\pdfoutput=1
\documentclass[11pt]{article}
\usepackage{authblk}
\usepackage[hyperfootnotes=false]{hyperref}
\usepackage{algorithm}
\usepackage{algpseudocode}
\usepackage{notation-2}
\usepackage{arxiv-2}
\usepackage{amssymb}
\usepackage{enumitem}
\usepackage[T1]{fontenc}
\usepackage{lmodern}
\usepackage{microtype}
\hypersetup{colorlinks=true,allcolors=blue}
\usepackage{float}

\title{%
  \normalfont
  \makebox[\textwidth][c]{%
    \parbox{1\textwidth}{%
      \centering
      Beyond Second-Order Stationarity: Finding Local Minimax Points with Gradient and Comparison Oracles
    }%
  }%
}

\author[1,2]{Haichen Hu\thanks{\texttt{huhc@mit.edu}}}
\author[2,3]{David Simchi-Levi\thanks{\texttt{dslevi@mit.edu}}}
\affil[1]{Center for Computational Science and Engineering, MIT}
\affil[2]{Department of Civil and Environmental Engineering, MIT}
\affil[3]{Institute for Data, Systems, and Society, MIT}
\date{}

\begin{document}
\maketitle
\begin{abstract}
We study first-order methods for finding local minimax
optima of smooth nonconvex-nonconcave objectives.
We provide a partial resolution of the open problem
posed by \citet{chae2023open} for nondegenerate calm
local minimax points \citep{ma2025calm}.
Calmness bounds the maximizer's response radius
proportionally to the minimizer's perturbation radius.
We develop the Comparison-Augmented Gradient Descent Ascent
(\texttt{CA-GDA}) algorithm, which combines first-order gradient
updates with zeroth-order function comparisons. Each comparison
can be implemented using two function-value evaluations, so the
algorithm is implementable with the standard first-order oracle.
We prove local convergence to every nondegenerate
calm local minimax point, including non-strict optima
with singular maximizer Hessians. Convergence is geometric
after finitely many restarts.
If every stationary point is nondegenerate and at least
one calm local minimax point exists, we further establish
almost-sure convergence of the entire iterate sequence
to a calm local minimax point from any nonstationary
initialization.
These results establish attraction to every point in the desired
solution class and, under the global assumptions, convergence
exclusively to that class. The guarantees therefore go beyond
second-order stationarity: they identify local minimax optima
themselves, rather than merely stationary points satisfying
second-order necessary conditions.
\end{abstract}

\section{Introduction}
\label{sec:introduction}

Minimax optimization, $\min_{x}\max_{y}f(x,y)$, provides a fundamental mathematical
framework for decision-making under uncertainty and
strategic interaction, with applications across game theory,
statistics, operations research, and control.
It describes competition between players with opposing
objectives, the design of statistical procedures that
minimize worst-case risk \citep{wald1950statistical},
decisions that remain effective under uncertain parameters
\citep{bental2009robust}, and controllers that maintain
performance against adverse disturbances
\citep{basar1995hinfinity}.
Across these applications, the minimization variable
represents a decision, while the maximization variable
represents an opposing strategy or an adverse scenario.
For nonconvex-nonconcave objectives, the interaction between
minimization and maximization complicates both the definition
of local optimality and the design of algorithms that attain it.

Local optimality in a minimax problem must account for
its sequential structure: the minimizer chooses an action,
and the maximizer responds to that choice.
To capture this structure, \citet{jin2020local} introduced
the notion of a local minimax point.
Specifically, a point $(x^*,y^*)$ is local minimax if
there exist $\delta_0>0$ and a response-radius function
$h:(0,\delta_0]\to(0,\infty)$ satisfying
$h(\delta)\to0$ as $\delta\downarrow0$, such that,
for every $\delta\in(0,\delta_0]$ and every $x,y$ with
$\norm{x-x^*}\le\delta$ and $\norm{y-y^*}\le\delta$,
\[
    f(x^*,y)
    \le f(x^*,y^*)
    \le \max_{\norm{y'-y^*}\le h(\delta)} f(x,y').
\]
The first inequality requires $y^*$ to locally maximize
the objective at $x^*$.
The second allows the maximizer to respond to a perturbation
of $x^*$ within a distance $h(\delta)$ of $y^*$ and to recover
at least the original objective value.
The requirement $h(\delta)\to0$ ensures that this response
remains local as the minimizer's perturbation vanishes.

The local stability guarantee for gradient descent
ascent in \citet{jin2020local} requires the target to satisfy
the strict local minimax conditions.
Specifically, at a stationary point $(x^*,y^*)$, these
conditions require $\nabla_{yy}^2 f \prec 0,\ \nabla_{xx}^2 f
    -\nabla_{xy}^2 f
     (\nabla_{yy}^2 f)^{-1}
     \nabla_{yx}^2 f \succ 0$, where all derivatives are evaluated at $(x^*,y^*)$.
The extragradient methods of \citet{chae2024two}
and \citet{kim2024double} extend convergence guarantees
to additional non-strict points.
However, in the general non-strict setting, these guarantees
identify stationary points satisfying refined second-order
necessary conditions for local minimax optimality.
They do not ensure that the limit satisfies the definition
of a local minimax point.
Indeed, the set of stationary points satisfying these
necessary conditions can be strictly larger than the set
of local minimax points: even with a nonsingular full Hessian,
such a point may fail to locally maximize the objective
in the maximization variable.
Thus, convergence to the set defined by these necessary
conditions leaves unresolved whether a first-order method
can converge exclusively to local minimax optima themselves.
This gap is central to the open problem posed by
\citet{chae2023open}:
\begin{quote}
    \itshape
    Is there a first-order method that only converges
    to local minimax optima?
\end{quote}

In this paper, we provide a partial resolution of this
open problem for nondegenerate calm local minimax points.
Calm local minimax optimality, introduced by
\citet{ma2025calm}, requires that the maximizer's response
radius be chosen proportional to the minimizer's
perturbation radius.
Equivalently, the response radius satisfies
$h(\delta)=O(\delta)$.
This restriction also appears in the optimality conditions
discussed by \citet{chae2023open}.
The class includes non-strict local minimax points:
the nonsingularity of the full Hessian
$\nabla^2 f(x^*,y^*)$ permits
$\nabla_{yy}^2 f(x^*,y^*)$ to be singular.

Our approach combines gradient updates, stabilized by an
auxiliary variable, with randomized function comparisons.
The gradient updates provide local attraction to the desired
points, while the comparisons detect improvements that violate
local maximality in the maximization variable.
An adaptive comparison margin allows persistent improvements
to be detected while preserving convergence near a solution.
A restart procedure selects update parameters and sampling
radii without knowing the solution. It first restarts near
the current iterate and, when needed, samples new
initializations throughout the space.
Our main contributions are as follows.

\begin{enumerate}
    \item
    We develop an algorithm using one gradient query and
    one single-bit function comparison per candidate update.
    Each comparison can be implemented with two function-value
    evaluations, so the method is implementable with the standard
    first-order oracle.

    \item
    We establish local convergence to every nondegenerate
    calm local minimax point
    (Theorem~\ref{thm:local-convergence}).
    From every sufficiently close initialization, the
    algorithm either terminates at that point
    or converges to it geometrically after finitely many
    restarts.
    Nondegeneracy is required only at the target point,
    whose maximizer Hessian may be singular.

    \item
    We establish almost-sure global convergence when every
    stationary point is nondegenerate and at least one calm
    local minimax point exists
    (Theorem~\ref{thm:global-convergence}).
    From any nonstationary initialization, the algorithm
    makes finitely many restarts almost surely, and its
    entire sequence converges to a calm local minimax point.
    Thus, the guarantee identifies actual calm local minimax
    points, rather than a larger set characterized only by
    necessary optimality conditions.
\end{enumerate}

Together, these results establish local attraction to
nondegenerate calm local minimax points and almost-sure
global convergence to this class under the respective
assumptions.
They provide a partial resolution of the open problem
posed by \citet{chae2023open}, covering non-strict points
with singular maximizer Hessians.

\paragraph{Notations.}
For an open set $\mathcal U$, $C^k(\mathcal U)$ denotes
the class of functions whose derivatives up to order $k$
exist and are continuous.
We write $\nabla f$ and $\nabla^2 f$ for the gradient
and Hessian of $f$, with subscripts indicating the
variables of differentiation.
For a differentiable map $G$, $DG$ denotes its Jacobian.
We use $\norm{\cdot}$ for the Euclidean norm of vectors
and the induced operator norm of matrices.
For a matrix $M$, $M^\top$, $\ker M$,
$\operatorname{range}M$, and $\operatorname{rank}M$
denote its transpose, kernel, range, and rank.
The identity matrix is denoted by $I_d$, or simply $I$
when its dimension is clear.
For symmetric matrices, $\succeq0$ and $\succ0$ denote
positive semidefiniteness and positive definiteness;
$\preceq0$ and $\prec0$ are defined analogously.
For a square matrix $M$, $\rho(M)$ denotes its spectral
radius.
For a nonempty set $S$, we define
$\operatorname{dist}(z,S)=\inf_{u\in S}\norm{z-u}$.
We write $\mathbf{1}\{E\}$ for the indicator of an event $E$ and
$\operatorname{Leb}$ for Lebesgue measure.
The notation $O(\cdot)$ and $o(\cdot)$ is standard Landau notation.
\section{Related Works}
\paragraph{Local minimax optimization.}
\citet{jin2020local} introduced local minimax points
and connected their optimality with the stability
of two-timescale gradient descent ascent.
Subsequent work developed optimality conditions,
compared local solution concepts, and introduced
calm local minimax optimality
\citep{zhang2022optimality,ma2025calm}.
Algorithmic approaches include implicit Stackelberg
learning, Hessian-based corrections, and finite
timescale separation
\citep{fiez2020implicit,wang2020solving,fiez2021local}.
More recently, extragradient methods extended stability
guarantees to non-strict points with singular maximizer
Hessians \citep{chae2024two,kim2024double}.
Their curvature conditions alone do not ensure local
maximality in the maximization variable. This gap is discussed
by \citet{chae2023open}.
Related global convergence results target min-max
critical points, a different solution concept
\citep{daskalakis2023stayon}.
Our work combines gradient updates with function
comparisons to establish convergence to nondegenerate
calm local minimax points themselves.

\paragraph{First-order optimization.}
Classical first-order frameworks include accelerated
gradient, mirror descent, and proximal-gradient methods
\citep{nesterov2004introductory,beck2003mirror,
parikh2014proximal}.
For nonconvex optimization, stochastic approximation
and variance reduction provide efficient methods for
computing approximate stationary points
\citep{ghadimi2013stochastic,reddi2016stochastic,
fang2018spider}.
Online-to-nonconvex conversions \citep{ahn2024general,cutkosky2023optimal,hu2026optimizing}
provide another approach to stochastic nonconvex optimization
by translating online learning guarantees into stationarity guarantees.
Oracle lower bounds characterize the worst-case
complexity of finding approximate stationary points
\citep{carmon2020lower}.
Random initialization and perturbations further support
strict-saddle avoidance and approximate second-order
stationarity under suitable regularity conditions
\citep{lee2016gradient,lee2019first,
panageas2019first,jin2017escape}.
For monotone variational inequalities and structured
saddle-point problems, mirror-prox, extragradient,
and optimistic updates provide convergence guarantees
through proximal and extrapolation mechanisms
\citep{nemirovski2004prox,daskalakis2018training,
mokhtari2020unified}.
General nonconvex-nonconcave games present additional
difficulties, including attraction to cycles and other
sets containing no stationary points
\citep{hsieh2021limits}.
Auxiliary-variable smoothing and dissipative
regularization can stabilize gradient descent ascent
under structural assumptions
\citep{zhang2020single,zheng2023universal,
zheng2024dissipative}.
Our auxiliary-variable gradient update is related to
these stabilization schemes.
We additionally use function comparisons to detect
violations of local maximality in the maximization
variable and establish convergence to nondegenerate
calm local minimax points.

\paragraph{Zeroth-order optimization and comparison oracles.}
Zeroth-order methods optimize objectives through
function evaluations, using randomized smoothing
and gradient estimation to construct search directions
\citep{ghadimi2013stochastic,duchi2015optimal,
nesterov2017random}.
Randomized direct-search methods also select iterates
by comparing the objective values of sampled candidates
\citep{bergou2020stochastic}.
Comparison-based optimization studies a more restricted
feedback model in which queries reveal relative
preferences instead of numerical function values.
Algorithms and query-complexity guarantees have been
developed for noisy pairwise comparisons, noiseless
comparator oracles, and dueling convex optimization
\citep{jamieson2012query,karabag2021smooth,
saha2021dueling,saha2025dueling}.
Related methods estimate normalized gradient directions
from one-bit comparison feedback and use these estimates
to guide descent \citep{cai2022onebit}.
More recently, \citet{scheinberg2026functionfree}
developed a function-free optimization framework based
directly on preference relations and level-set geometry.
Our algorithm combines gradient updates with comparisons
that test whether a proposed maximizer action improves
the objective by more than an adaptive additive margin.
We use this feedback to detect violations of local
maximality while preserving convergence near calm
local minimax points.
\section{Problem Setup}
\label{sec:model}

In this section, we introduce the optimization problem
and review the definitions and optimization background
used in the analysis.

Let $f:\RR^{d_1}\times\RR^{d_2}\to\RR$ be a function,
where $d_1$ and $d_2$ are positive integers. We consider the
following minimax optimization problem
\[
    \min_{x\in\RR^{d_1}}\max_{y\in\RR^{d_2}} f(x,y).
\]

We write $z=(x,y)$ and use $f(z)$ and $f(x,y)$ interchangeably.
Throughout the paper, $\norm{\cdot}$ denotes the Euclidean norm
for vectors and the induced operator norm for matrices.
Following standard smoothness assumptions in minimax optimization
\citep{chae2024two,kim2024double}, we make the following assumption on $f$.

\begin{assumption}
\label{ass:smooth}
The objective function $f$ belongs to $C^2(\RR^{d_1+d_2})$,
and there exists a constant $L>0$ such that
$\norm{\nabla^2 f(z)}\le L$ for every
$z\in\RR^{d_1+d_2}$.
\end{assumption}

We distinguish two types of stationary points of $f$.
\begin{definition}\label{def:stationary}
    For any $z_0$ such that $\nabla_z f(z)|_{z=z_0}=0$, if $\det(\nabla^2_zf(z))|_{z=z_0}\neq 0$, then we say that $z_0$ is a nondegenerate stationary point. Otherwise, we call $z_0$ a degenerate stationary point.
\end{definition}
We focus on nondegenerate local minimax optimal points.

We next specify the notion of local optimality considered in this
paper. In the minimax problem, the minimizer first chooses an
action $x$, and the maximizer responds with an action $y$.
A local solution should therefore protect the value against small
changes in $y$ at the proposed $x$, and allow the maximizer to
recover that value when the minimizer moves.
We focus on \emph{calm local minimax points}, introduced by
\citet{ma2025calm}.

\begin{definition}[\cite{ma2025calm}]
\label{def:calm_local_minimax}
A point $z^*=(x^*,y^*)$ is a \emph{calm local minimax point}
if there exist constants $M>0$ and $\delta_0>0$
and a radius function $h(\delta)=M\delta$ such that,
for every $\delta\in(0,\delta_0]$ and every $x,y$ satisfying
$\norm{x-x^*}\le\delta$ and $\norm{y-y^*}\le\delta$,
\begin{equation}
    f(x^*,y)\le f(x^*,y^*)\le
    \max_{\norm{y'-y^*}\le h(\delta)} f(x,y').
    \label{eq:calm-definition}
\end{equation}
\end{definition}

The constant $M$ bounds the size of the maximizer's response
neighborhood relative to the minimizer's perturbation radius,
uniformly over all sufficiently small $\delta$.
Equivalently, we may use a valid response radius
$h(\delta)=O(\delta)$, since enlarging the response ball can only
increase the maximum.
If we only require $h(\delta)\to0$ as $\delta\downarrow0$,
we recover the original local minimax definition of
\citet{jin2020local}.

Calm local minimax points include non-strict, nondegenerate
points with singular $\nabla_{yy}^2 f(z^*)$.
For example, consider $f(x,y)=-x^2+2xy$.
At the origin, $f(0,y)=0$ for every $y$.
For any $|x|\le\delta$, choosing $y'=x/2$ gives
$f(x,y')=0$ and $|y'|\le\delta/2$.
Thus, the origin is a calm local minimax point with $M=1/2$,
although $\nabla_{yy}^2 f(0,0)=0$.
It is also nondegenerate since $\nabla^2 f(0,0)=
    \begin{pmatrix}
        -2 & 2\\
         2 & 0
    \end{pmatrix}$.

The open problem posed by \citet{chae2023open} asks whether
there exists a first-order method that converges only to local
minimax optima. We therefore clarify the meaning of
zeroth-order, first-order, and second-order methods by specifying
the information available through their oracles
\citep{nemirovski1983problem,carmon2020lower}.

For a query point $z\in\RR^{d_1+d_2}$, a zeroth-order oracle
returns the function value, a first-order oracle additionally
returns the gradient, and a second-order oracle additionally
returns the Hessian. We write these oracles as
\begin{align*}
    \mathcal O_f^{(0)}(z)
    &= f(z),\\
    \mathcal O_f^{(1)}(z)
    &= \bigl(f(z),\nabla f(z)\bigr),\\
    \mathcal O_f^{(2)}(z)
    &= \bigl(f(z),\nabla f(z),\nabla^2 f(z)\bigr).
\end{align*}
A method that accesses $f$ through $\mathcal O_f^{(p)}$,
for $p\in\{0,1,2\}$, is called a zeroth-order, first-order,
or second-order method, respectively.
Its query points and iterates may depend on previous oracle
answers and its internal randomness.

Comparison oracles, which report the relative preference between
two queried points, have also been studied in derivative-free
optimization \citep{scheinberg2026functionfree}.
Our algorithm combines first-order information from a gradient
oracle with zeroth-order information from a function comparison
oracle with a margin.

\begin{definition}[Gradient oracle]
\label{def:gradient-oracle}
For any point $z=(x,y)\in\RR^{d_1+d_2}$,
the gradient oracle returns
\[
    \mathcal G_f(z)
    =\nabla f(z)
    =\begin{pmatrix}
        \nabla_x f(x,y)\\
        \nabla_y f(x,y)
    \end{pmatrix}.
\]
\end{definition}

The comparison oracle tests whether a proposed change in the
maximization variable improves the objective by more than a
specified margin.

\begin{definition}[Function comparison oracle with a margin]
\label{def:comparison-oracle}
Given $x\in\RR^{d_1}$, two points $y,y'\in\RR^{d_2}$,
and a margin $\vartheta\ge0$, the function comparison oracle
returns the bit
\begin{equation}
    \mathcal C_f(x;y,y',\vartheta)
    =\mathbf{1}\bigl\{
        f(x,y')>f(x,y)+\vartheta
    \bigr\}.
    \label{eq:comparison-oracle}
\end{equation}
The margin $\vartheta$ is chosen by the algorithm and may vary
across queries.
\end{definition}

Each iteration of our algorithm makes one gradient query and
one comparison query, with the latter returning a single bit.
A comparison query can be implemented using two
function-value evaluations. Our algorithm can therefore be
implemented using the standard first-order oracle.


\section{Comparison-Augmented Gradient Descent Ascent Algorithm}
\label{sec:algorithm}

In this section, we develop an algorithm that combines gradient
information with function comparisons to find nondegenerate calm
local minimax points. The algorithm first uses the gradient to
compute a candidate iterate. Holding the minimizer's proposed
action fixed, it then samples a nearby action for the maximizer
and uses the comparison oracle to test whether this action
increases the objective value by more than an adaptively
chosen margin.

To compute the candidate iterate, we use a modified gradient
descent-ascent update. Ordinary simultaneous gradient
descent-ascent can fail to converge even near a calm local
minimax point. For example, for the bilinear objective
$f(x,y)=xy$, the origin is a calm local minimax point,
but simultaneous gradient descent-ascent is unstable there
for every pair of positive constant step sizes. We therefore
introduce an auxiliary variable $v_t\in\RR^{d_2}$ that tracks
the maximizer's action $y_t$. The maximizer follows the
ascent direction $\nabla_y f(x_t,y_t)$ and receives
a correction $-s(y_t-v_t)$, where $s>0$.
This correction pulls its action toward $v_t$. The auxiliary variable gradually moves toward
$y_t$, while both $x_t$ and $v_t$ evolve on a slower timescale
than $y_t$. At a fixed point of this gradient update, $v=y$,
so the correction vanishes and the fixed-point equations
reduce to $\nabla f(x,y)=0$.

Specifically, write $w_t=(x_t,y_t,v_t)$ and set
$\eta=1/(4L)$.
For fixed parameters $s\in(0,L]$ and $\tau\ge1$, we query
$\nabla f(x_t,y_t)$ and compute
\begin{equation}
\begin{aligned}
    x_{t+1}^{\mathrm{cand}}
    &=x_t-\frac{\eta}{\tau}\nabla_x f(x_t,y_t),\\
    \widetilde y_{t+1}
    &=y_t+\eta\bigl(
        \nabla_y f(x_t,y_t)-s(y_t-v_t)
    \bigr),\\
    v_{t+1}^{\mathrm{cand}}
    &=v_t+\frac{\eta s}{\tau}(y_t-v_t).
\end{aligned}
\label{eq:base-update}
\end{equation}
The factor $1/\tau$ slows the updates of the minimizer
and the auxiliary variable relative to the maximizer's update
when $\tau>1$.
The superscript $\mathrm{cand}$ distinguishes a proposed
update from an accepted iterate.
The action $\widetilde y_{t+1}$ is provisional;
the comparison step below determines the maximizer's
proposed action $y_{t+1}^{\mathrm{cand}}$.

To construct an alternative to $\widetilde y_{t+1}$,
fix a sampling radius $\ell>0$ and draw $U_t$ independently
and uniformly from the unit ball in $\RR^{d_2}$.
The candidate action is
\begin{equation}
    \widehat y_{t+1}
    =\widetilde y_{t+1}+\ell U_t.
    \label{eq:probe-candidate}
\end{equation}
Using the gradient already obtained at the current iterate,
we set the comparison margin to
\begin{equation}
    \vartheta_t
    =
    \ell\sqrt{
        L\ell\left(
            \norm{\nabla f(x_t,y_t)}
            +s\norm{y_t-v_t}
        \right)
    }.
    \label{eq:adaptive-margin}
\end{equation}
We then hold $x_{t+1}^{\mathrm{cand}}$ fixed and query
the comparison oracle:
\begin{equation}
    y_{t+1}^{\mathrm{cand}}
    =
    \begin{cases}
        \widehat y_{t+1},
        & \text{if }
        \mathcal C_f\bigl(
            x_{t+1}^{\mathrm{cand}};
            \widetilde y_{t+1},
            \widehat y_{t+1},
            \vartheta_t
        \bigr)=1,\\
        \widetilde y_{t+1},
        & \text{otherwise}.
    \end{cases}
    \label{eq:comparison-update}
\end{equation}
Thus, the maximizer adopts the candidate action only when
\[
    f(x_{t+1}^{\mathrm{cand}},\widehat y_{t+1})
    >
    f(x_{t+1}^{\mathrm{cand}},\widetilde y_{t+1})
    +\vartheta_t.
\]
The comparison changes only the maximizer's proposed action.
The resulting candidate state is
\[
    w_{t+1}^{\mathrm{cand}}
    =
    \bigl(
        x_{t+1}^{\mathrm{cand}},
        y_{t+1}^{\mathrm{cand}},
        v_{t+1}^{\mathrm{cand}}
    \bigr).
\]

Intuitively, the gradient update uses information at the
current iterate. The comparison tests whether the maximizer
can obtain a larger value at a nearby point.
The margin requires this improvement to exceed
$\vartheta_t$, so an arbitrarily small increase does not
automatically trigger a move.
For fixed parameters, the margin tends to zero as
$\norm{\nabla f(x_t,y_t)}$ and $\norm{y_t-v_t}$ tend to zero.
A persistent positive improvement can therefore eventually
exceed the margin.
Computing the margin requires no additional gradient query.

The updates above use three parameters:
the correction strength $s$, the timescale parameter $\tau$,
and the sampling radius $\ell$.
Suitable values may depend on the unknown target.
In particular, local maximality of $y^*$ controls
$f(x^*,y)$ only within a neighborhood of $y^*$,
whose radius is not known to the algorithm.
We therefore select the parameters through a sequence of
trials. For a scale $\gamma\in(0,1)$, set
\begin{equation}
    s=L\gamma,\ 
    \tau=\gamma^{-2},\ 
    \ell=\gamma.
    \label{eq:parameter-scale}
\end{equation}
A smaller $\gamma$ reduces the sampling radius and the
correction strength while increasing the timescale separation.
All three parameters remain fixed within a trial.

We first attempt to find a solution by restarting near
the current iterate.
If these attempts repeatedly fail, the algorithm switches
to initialization centers sampled throughout the space.
Let $t_j$ denote the step
at which trial $j$ begins.
The index $t$ counts state updates across all trials,
including initialization perturbations.
Given $(x_0,y_0)$, set
\[
    w_0=(x_0,y_0,y_0),\  t_1=1.
\]
Compute
\begin{equation}
    \alpha_{\mathrm{init}}
    =
    \min\left\{
        \frac12,
        \left(
            \frac{\norm{\nabla f(x_0,y_0)}}{L}
        \right)^{1/3}
    \right\}.
    \label{eq:initial-restart-scale}
\end{equation}
If $\alpha_{\mathrm{init}}=0$, the method returns the
stationary point $(x_0,y_0)$.
Otherwise, define
\begin{equation}
    K=\left\lceil\alpha_{\mathrm{init}}^{-1}\right\rceil.
    \label{eq:local-trial-limit}
\end{equation}
The number $K$ is computed once and remains fixed.
It determines how many trials may restart near the
current iterate before the algorithm switches to
global exploration.

For each trial $j\le K$, use the retained state as the
initialization center and set
\begin{equation}
    \bar w_j=w_{t_j-1},\ 
    \gamma_j=2^{-j}.
    \label{eq:local-trial-center}
\end{equation}
If all $K$ trials fail, then, for each subsequent trial
$j>K$, independently draw
\begin{equation}
\begin{aligned}
    \bar w_j
    &\sim\mathcal N(w_0,I_{d_1+2d_2}),\\
    \Pr(\nu_j=m)
    &=2^{-m},\  m=1,2,\ldots,
\end{aligned}
\label{eq:global-trial-sampling}
\end{equation}
and set $\gamma_j=2^{-\nu_j}$.
The Gaussian distribution remains centered at the original
state $w_0$.
Conditional on $w_0$, the new center and parameter scale
are independent of each other and of all previous draws.
Sampling the scale anew gives each fixed value $2^{-m}$
repeated opportunities to be selected if trials continue.

In either case, write
$\bar w_j=(\bar x_j,\bar y_j,\bar v_j)$ and use
\eqref{eq:parameter-scale} with $\gamma=\gamma_j$
to obtain $s_j,\tau_j,\ell_j$.
At the selected center, compute
\begin{equation}
\begin{aligned}
    R_j
    &=
    \norm{\nabla f(\bar x_j,\bar y_j)}
    +s_j\norm{\bar y_j-\bar v_j},\\
    \alpha_j
    &=
    \min\left\{
        \frac12,
        \left(\frac{R_j}{L}\right)^{1/3}
    \right\}.
\end{aligned}
\label{eq:restart-scale}
\end{equation}
If $R_j=0$, the method returns the stationary point
$(\bar x_j,\bar y_j)$.
Otherwise, draw $Z_j$ independently and uniformly from
the unit ball in $\RR^{d_1+2d_2}$ and initialize the trial by
\begin{equation}
    w_{t_j}=\bar w_j+\alpha_j^3Z_j.
    \label{eq:trial-initialization}
\end{equation}

Starting at $t=t_j$, compute a candidate state
$w_{t+1}^{\mathrm{cand}}$ using
\eqref{eq:base-update} and \eqref{eq:comparison-update}
with parameters $s_j,\tau_j,\ell_j$.
While
\begin{equation}
    \norm{w_{t+1}^{\mathrm{cand}}-w_t}
    \le
    \alpha_j^2(1-\alpha_j)^{t-t_j},
    \label{eq:restart-check}
\end{equation}
accept the update by setting $    w_{t+1}=w_{t+1}^{\mathrm{cand}}$ and increase $t$ by one.
Then compute the next candidate.
The parameters and $\alpha_j$ remain fixed throughout
this trial.

If \eqref{eq:restart-check} first fails at step $t$,
discard the candidate and set $ t_{j+1}=t+1$.

Trial $j+1$ then selects its center and parameters
according to \eqref{eq:local-trial-center} or
\eqref{eq:global-trial-sampling}, recomputes
$\alpha_{j+1}$, and initializes
\[
    w_{t+1}
    =
    \bar w_{j+1}+\alpha_{j+1}^3Z_{j+1}.
\]
If the check continues to hold, the current trial
continues without a restart.
A trial that continues forever has finite total movement
under this bound, since
\[
    \sum_{n=0}^{\infty}
    \alpha_j^2(1-\alpha_j)^n
    =\alpha_j.
\]

Algorithm~\ref{alg:hybrid} summarizes the complete procedure
with a budget of $N$ state updates.
The convergence analysis considers the same procedure
continued without this finite budget.

\begin{algorithm}[htbp]
\small
\caption{Comparison-Augmented Gradient Descent Ascent (\texttt{CA-GDA})}
\label{alg:hybrid}
\begin{algorithmic}[1]
    \Require Initial point $(x_0,y_0)$, Hessian bound $L>0$,
    update budget $N\ge1$.
    \Statex \textit{All random samples are independent
    conditional on $w_0$.}

    \State $w_0\gets(x_0,y_0,y_0)$,
    $t\gets0$, $j\gets1$, $\eta\gets1/(4L)$
    \State $(g_x,g_y)\gets\mathcal G_f(x_0,y_0)$

    \If{$\norm{(g_x,g_y)}=0$}
        \State \Return $(x_0,y_0)$
    \EndIf

    \State $\alpha_{\mathrm{init}}\gets
    \min\{1/2,(\norm{(g_x,g_y)}/L)^{1/3}\}$,
    $K\gets\lceil1/\alpha_{\mathrm{init}}\rceil$

    \While{$t<N$}
        \If{$j\le K$}
            \State $(\bar x_j,\bar y_j,\bar v_j)\gets w_t$,
            $\gamma_j\gets2^{-j}$
        \Else
            \State Draw $(\bar x_j,\bar y_j,\bar v_j)$
            from $\mathcal N(w_0,I_{d_1+2d_2})$
            \State Draw $\nu_j$ with
            $\Pr(\nu_j=m)=2^{-m}$ for $m=1,2,\ldots$
            \State $\gamma_j\gets2^{-\nu_j}$
            \State $(g_x,g_y)\gets
            \mathcal G_f(\bar x_j,\bar y_j)$
        \EndIf

        \State $s_j\gets L\gamma_j$,
        $\tau_j\gets\gamma_j^{-2}$,
        $\ell_j\gets\gamma_j$
        \State $R_j\gets
        \norm{(g_x,g_y)}
        +s_j\norm{\bar y_j-\bar v_j}$

        \If{$R_j=0$}
            \State \Return $(\bar x_j,\bar y_j)$
        \EndIf

        \State $\alpha_j\gets
        \min\{1/2,(R_j/L)^{1/3}\}$
        \State $t_j\gets t+1$
        \State Draw $Z_j$ uniformly from the unit ball
        in $\RR^{d_1+2d_2}$
        \State $w_{t_j}\gets
        (\bar x_j,\bar y_j,\bar v_j)+\alpha_j^3Z_j$
        \State $t\gets t_j$

        \While{$t<N$}
            \State $(g_x,g_y)\gets\mathcal G_f(x_t,y_t)$
            \State $a_t\gets s_j(y_t-v_t)$

            \State
            $(x_{t+1}^{\mathrm{cand}},v_{t+1}^{\mathrm{cand}})
            \gets
            (x_t-\eta g_x/\tau_j,\,
             v_t+\eta a_t/\tau_j)$
            \State $\widetilde y_{t+1}
            \gets y_t+\eta(g_y-a_t)$

            \State $\vartheta_t\gets
            \ell_j\sqrt{
                L\ell_j\bigl(
                    \norm{(g_x,g_y)}+\norm{a_t}
                \bigr)
            }$
            \State Draw $U_t$ uniformly from the unit ball
            in $\RR^{d_2}$
            \State $\widehat y_{t+1}
            \gets\widetilde y_{t+1}+\ell_jU_t$

            \State $b_t\gets
            \mathcal C_f\bigl(
                x_{t+1}^{\mathrm{cand}};
                \widetilde y_{t+1},
                \widehat y_{t+1},
                \vartheta_t
            \bigr)$
            \State $y_{t+1}^{\mathrm{cand}}
            \gets\widetilde y_{t+1}+b_t\ell_jU_t$
            \State $w_{t+1}^{\mathrm{cand}}
            \gets\bigl(
                x_{t+1}^{\mathrm{cand}},
                y_{t+1}^{\mathrm{cand}},
                v_{t+1}^{\mathrm{cand}}
            \bigr)$

            \If{$\norm{w_{t+1}^{\mathrm{cand}}-w_t}
                \le\alpha_j^2(1-\alpha_j)^{t-t_j}$}
                \State $w_{t+1}\gets w_{t+1}^{\mathrm{cand}}$,
                $t\gets t+1$
            \Else
                \State \textbf{break}
                \Comment{Discard the candidate and retain $w_t$}
            \EndIf
        \EndWhile

        \State $j\gets j+1$
    \EndWhile

    \State \Return $(x_t,y_t)$
\end{algorithmic}
\end{algorithm}

Each candidate update uses one gradient query and one
comparison query.
The gradient is reused to compute the comparison margin.
When a trial fails, its last gradient query is also reused
to initialize the next trial if $j+1\le K$.
In addition, the method makes one initial gradient query
and one gradient query at each newly sampled center
during global exploration.

\section{Theoretical Guarantees}
\label{sec:theory}
In this section, we provide theoretical guarantees for our
algorithm \texttt{CA-GDA}. Specifically, we first establish local convergence
to a nondegenerate calm local minimax point when initialized
sufficiently close to it. We then establish almost-sure global
convergence from any nonstationary initial point, under the
additional assumptions that all stationary points are
nondegenerate and at least one calm local minimax point exists.

We begin with the local convergence result.
Let $z^*=(x^*,y^*)$ be a calm local minimax point such that
$\nabla^2 f(z^*)$ is nonsingular, and write
$w^*=(x^*,y^*,y^*)$ for the corresponding augmented state. The nondegeneracy condition at $z^*$ is also used in
existing local stability analyses of minimax methods
\citep{chae2024two,kim2024double}.
This condition concerns the full Hessian
and allows $\nabla_{yy}^2 f(z^*)$ to be singular.
No nondegeneracy assumption is imposed on other stationary
points.
The following theorem shows that, when initialized
sufficiently close to $z^*$, the algorithm converges
to $w^*$ after finitely many restarts,
without entering the global exploration stage. For the convergence analysis,
we run Algorithm~\ref{alg:hybrid} without a prescribed iteration limit.

\begin{theorem}[Local convergence]
\label{thm:local-convergence}
Suppose Assumption~\ref{ass:smooth} holds.
Let $z^*=(x^*,y^*)$ be a calm local minimax point such that
$\nabla^2 f(z^*)$ is nonsingular, and write
$w^*=(x^*,y^*,y^*)$.

There exist $r>0$ and an integer $J_*\ge1$ with the
following property.
Initialize Algorithm~\ref{alg:hybrid} at
$w_0=(x_0,y_0,y_0)$, where
\[
    \norm{(x_0,y_0)-z^*}<r.
\]
The algorithm makes at most $J_*-1$ restarts.
Every trial uses the initialization rule
in \eqref{eq:local-trial-center}.

The algorithm either terminates at $z^*$ or converges to $w^*$.
If it does not terminate, let $J\le J_*$ denote its
final trial. Then
\begin{equation}
    \norm{w_{t_J+n}-w^*}
    \le
    \alpha_J(1-\alpha_J)^n,
    \  n=0,1,2,\ldots.
    \label{eq:local-convergence-rate}
\end{equation}
All conclusions hold for every choice of $(Z_j)_{j\ge1}$
and $(U_t)_{t\ge0}$ satisfying
$\norm{Z_j}\le1$ and $\norm{U_t}\le1$.
\end{theorem}
We next consider convergence from initial points that need
not lie near a solution.
For this result, we assume that at least one calm local
minimax point exists and that every stationary point of $f$
is nondegenerate.

\begin{theorem}[Global convergence]
\label{thm:global-convergence}
Suppose Assumption~\ref{ass:smooth} holds.
Assume that $f$ admits at least one calm local minimax point
and that every stationary point of $f$ is nondegenerate.

For any initial point $(x_0,y_0)$ satisfying
$\nabla f(x_0,y_0)\ne0$, initialize
Algorithm~\ref{alg:hybrid} at $w_0=(x_0,y_0,y_0)$.
Then, with probability one, the algorithm makes finitely
many restarts, and there exists a calm local minimax point
$z^\infty=(x^\infty,y^\infty)$ such that
\[
    w_t\longrightarrow
    w^\infty=(x^\infty,y^\infty,y^\infty).
\]
Moreover, if $J$ denotes the final trial, then
\begin{equation}
    \norm{w_{t_J+n}-w^\infty}
    \le
    \alpha_J(1-\alpha_J)^n,
    \  n=0,1,2,\ldots.
    \label{eq:global-convergence-rate}
\end{equation}
\end{theorem}
Both theorems establish geometric convergence in the
final trial, reached after finitely many restarts.

Theorem~\ref{thm:global-convergence} guarantees almost-sure
convergence of the entire sequence to a calm local minimax
point. The second-order necessary conditions studied
in \citet{chae2024two,kim2024double} may also hold at points
that are not local minimax optima. Our guarantee excludes
such limits.

Together with Theorem~\ref{thm:local-convergence}, our results
provide both directions under the assumptions of
Theorem~\ref{thm:global-convergence}.
Every calm local minimax point attracts sufficiently close
initializations for every allowed choice of the random vectors
$Z_j$ and $U_t$. From any nonstationary initial point,
the algorithm almost surely converges to a point in this class.
Thus, the guarantees identify the desired local optima
themselves, rather than a larger set of stationary points
defined only by second-order necessary conditions.

\section{Numerical Experiments}
\label{sec:numerical_experiment}

In this section, we evaluate our algorithm through three
sets of numerical experiments.
The adaptive parameter selection and restart procedure in
Algorithm~\ref{alg:hybrid} allow us to establish convergence
without knowing suitable parameter values in advance.
This includes almost-sure global convergence under the stated
assumptions.
For simplicity, our experiments directly use the gradient
and comparison updates with fixed step size $\eta$,
correction strength $s$, timescale parameter $\tau$,
and sampling radius $\ell$, without restarts.
The comparison margin remains adaptive.

We first compare our algorithm with two existing
extragradient methods.
We then remove the comparison step to examine the role
of the comparison oracle.
Finally, we compare the adaptive margin with a zero margin
to examine how the margin affects convergence near
a valid solution.
All convergence plots use the iteration number on the
horizontal axis and the Euclidean distance to the
corresponding calm local minimax solution set on the
vertical axis, displayed on a logarithmic scale.
Our curves report medians over 30 runs.

\paragraph{Comparison with extragradient methods.}
We compare our algorithm with two-timescale extragradient
($\tau$-EG) \citep{chae2024two} and double-step alternating
extragradient with increasing timescale separation
(Alt2-EG-ITS) \citep{kim2024double}.
We consider the scalar objective
\begin{equation}
    f_1(x,y)=\frac{x^2}{2}+xy+p_\rho(y),
    \  x,y\in\RR,
    \label{eq:numerical-first-function}
\end{equation}
and the objective
\begin{equation}
\begin{aligned}
    f_2(x,y)={}&-\frac{x_1^2}{2}
    +x_1y_1-p_\rho(y_1)+\frac{x_2^2}{2}
    +x_2y_2+p_\rho(y_2),
    \  x,y\in\RR^2,
\end{aligned}
\label{eq:numerical-second-function}
\end{equation}
where
$p_\rho(u)=2\rho^2(1-\cos(u/\rho))^2$
and $\rho=0.1$.
Let $a=\rho r_+\approx0.277464$, where $r_+$ is the larger
positive solution of $4(1-\cos r)\sin r=r$.
The corresponding calm local minimax solution sets are
\[
\begin{aligned}
    \mathcal Z_1
    &=\{(-a,a),(a,-a)\},\\
    \mathcal Z_2
    &=\bigl\{
        ((0,-\sigma a),(0,\sigma a)):
        \sigma\in\{-1,1\}
    \bigr\}.
\end{aligned}
\]
All these solutions are nondegenerate.
The solutions of $f_2$ have singular maximizer
Hessians and are therefore non-strict.
For both objectives, the origin is also a nondegenerate
stationary point, but it fails local maximality in $y$.
Indeed, $f_1(0,y)=p_\rho(y)>0$ and
$f_2(0,(0,u))=p_\rho(u)>0$ for sufficiently small nonzero
$y$ and $u$.
We use these objectives to test whether
an algorithm can avoid an incorrect stationary point
and reach a calm local minimax solution.

We run each method for $300{,}000$ iterations from
30 paired initializations near the origin.
Within each run, all three methods start from the same
$(x_0,y_0)$, sampled uniformly from the sphere
\[
    \norm{(x_0,y_0)}
    =10^{-4}\sqrt{d_1+d_2}.
\]
For our algorithm, we initialize $v_0=y_0$ and use
$L=9$, $\eta=1/36$, $s=9/16$, $\tau=256$,
and $\ell=1/16$.
The $\tau$-EG baseline uses $\eta=1/36$ with $\tau=20$
on $f_1$ and $\tau=200$ on $f_2$.
Alt2-EG-ITS uses its prescribed schedule
$\tau_t=(t+1)^{1/2.2}$ and alternating scale
$\sqrt{\tau_t}$, with step size $1/10$ on $f_1$
and $1/36$ on $f_2$.

\begin{figure}[h]
    \centering
    \includegraphics[width=\linewidth]{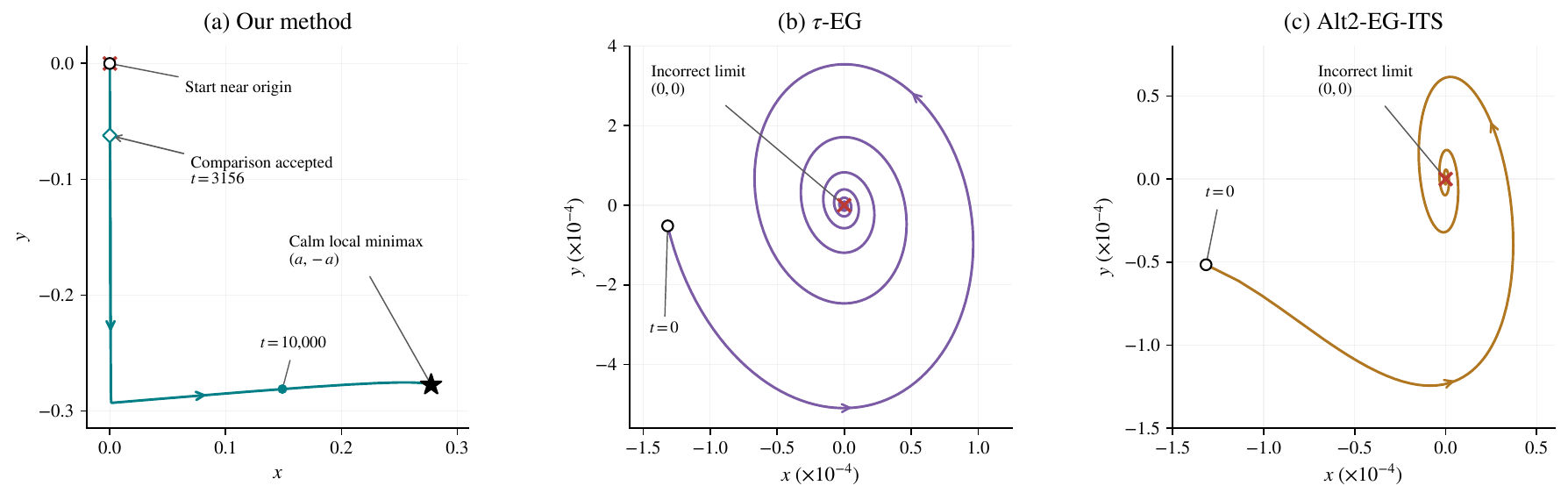}
    \caption{Trajectories for $f_1$:
    our algorithm, $\tau$-EG, and Alt2-EG-ITS
    from left to right.}
    \label{fig:numerical-first}
\end{figure}

Figure~\ref{fig:numerical-first} displays representative
trajectories for $f_1$ in the $(x,y)$-plane.
From left to right, the panels show our algorithm,
$\tau$-EG, and Alt2-EG-ITS.
Our algorithm leaves the neighborhood of the origin
and approaches the calm local minimax point $(a,-a)$.
Both extragradient methods instead approach the
incorrect stationary point $(0,0)$.
Arrows indicate the direction of the updates.
The panels use different coordinate ranges so that
the trajectories approaching the origin remain visible.

\begin{figure}[h]
    \centering
    \includegraphics[width=\linewidth]{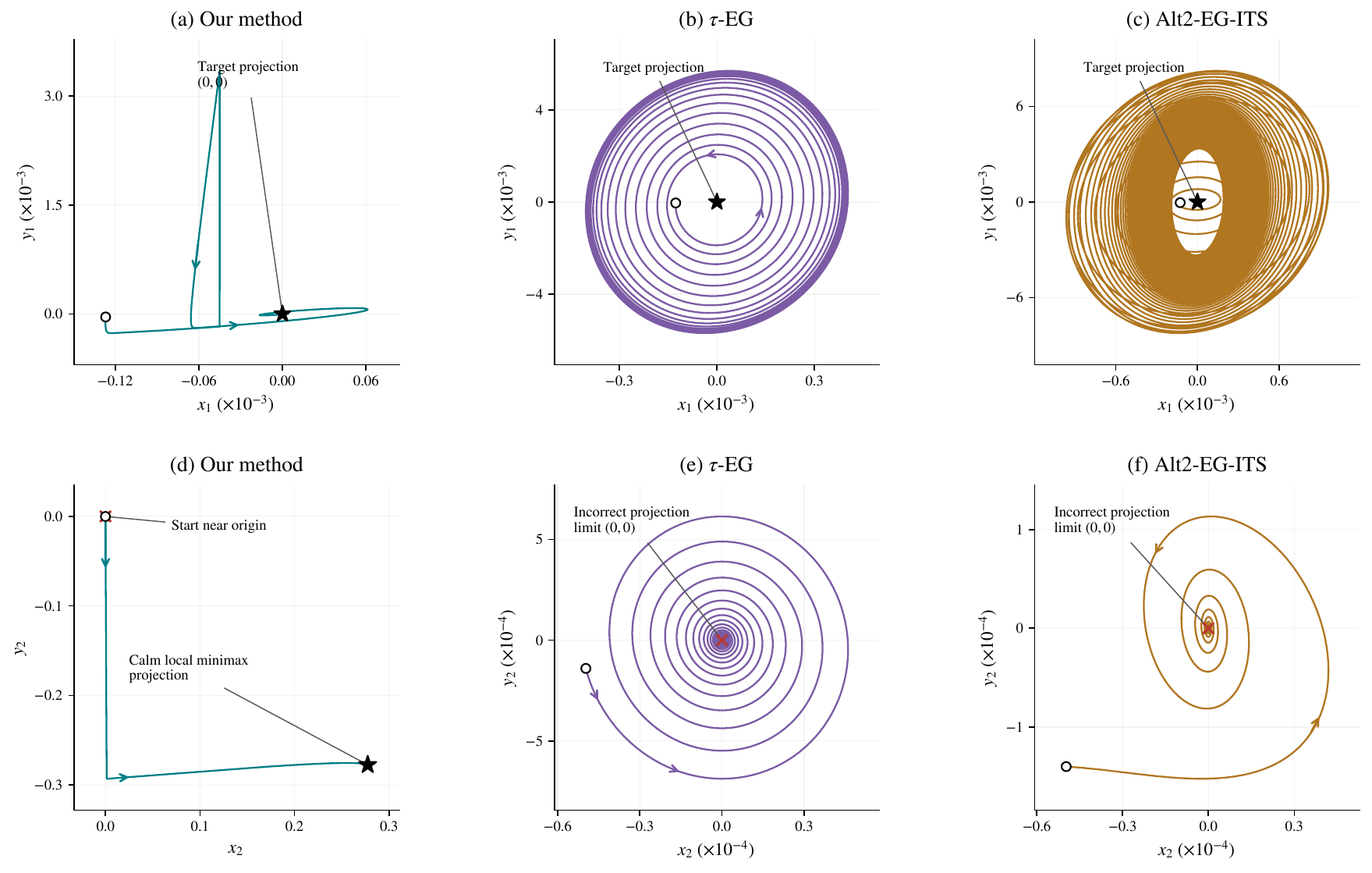}
    \caption{Projected trajectories for $f_2$:
    $(x_1,y_1)$ in the upper row and
    $(x_2,y_2)$ in the lower row.}
    \label{fig:numerical-second}
\end{figure}

Figure~\ref{fig:numerical-second} displays representative
trajectories for $f_2$.
The upper row shows the $(x_1,y_1)$ projections, and
the lower row shows the $(x_2,y_2)$ projections.
The columns follow the same method order as
Figure~\ref{fig:numerical-first}, with separate coordinate
ranges for each panel.
A calm local minimax solution requires the first
coordinate pair to approach $(0,0)$ and the second
to approach either $(-a,a)$ or $(a,-a)$.
Our algorithm exhibits both behaviors.
For the extragradient methods, the first coordinate
pair oscillates, while the second approaches the
incorrect projection $(0,0)$.
Thus, neither baseline reaches the calm local minimax
solution set in these runs.

\begin{figure}[h]
    \centering
    \includegraphics[width=0.8\linewidth]{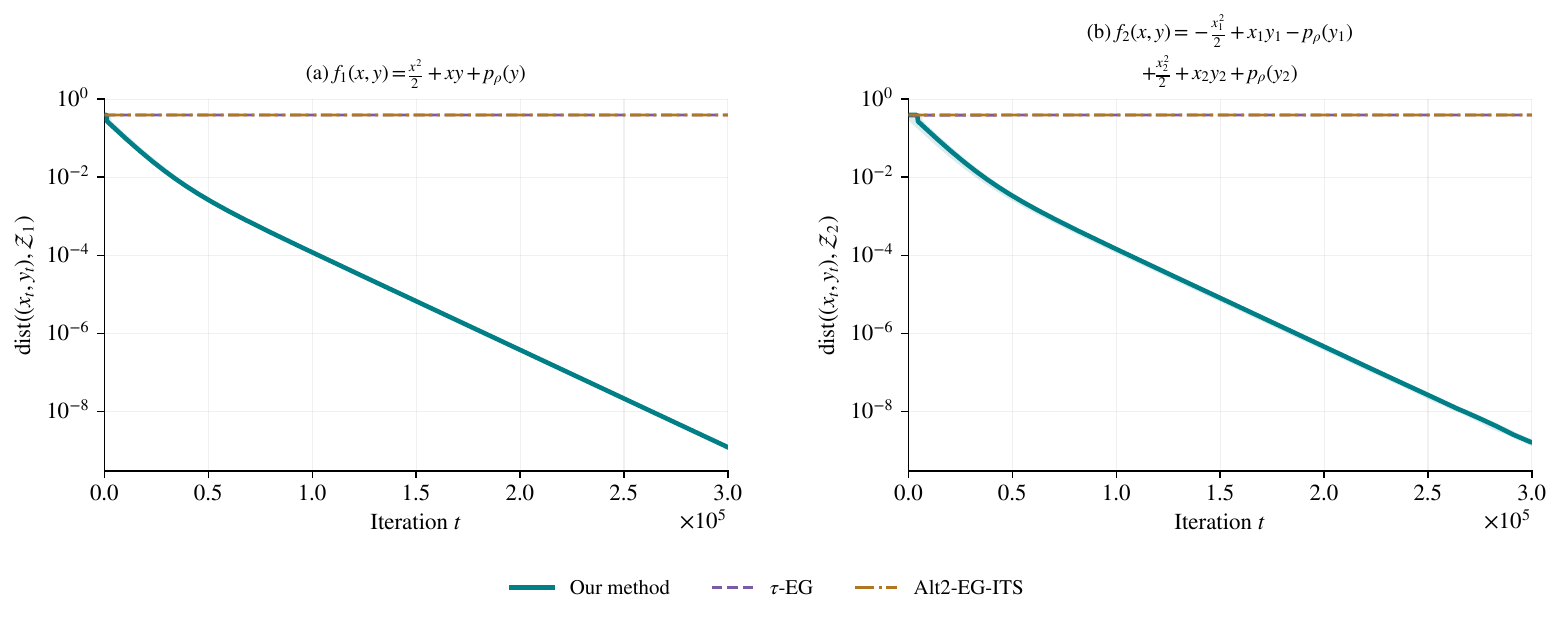}
    \caption{Convergence to $\mathcal Z_1$ for $f_1$
    (left) and $\mathcal Z_2$ for $f_2$ (right).}
    \label{fig:benchmark-convergence}
    \label{fig:numerical-first-convergence}
    \label{fig:numerical-second-convergence}
\end{figure}

Figure~\ref{fig:benchmark-convergence} summarizes
convergence over all 30 runs.
The left panel plots
$\operatorname{dist}((x_t,y_t),\mathcal Z_1)$
for $f_1$, and the right panel plots
$\operatorname{dist}((x_t,y_t),\mathcal Z_2)$
for $f_2$.
For $f_2$, this distance uses the complete
four-dimensional iterate rather than either projection.
In these runs, our algorithm approaches the local minimax
solution sets, while the two extragradient methods remain
separated from them.

On the tested problems and initializations, our algorithm
approaches calm local minimax solutions, including
non-strict solutions of $f_2$, while both extragradient
baselines remain separated from the solution sets.

\paragraph{Role of the comparison step.}
We next examine whether the comparison oracle is
responsible for avoiding the incorrect stationary points
in $f_1$ and $f_2$.
We remove only the comparison step, retaining the same
gradient updates for $x$, $y$, and the auxiliary variable
$v$, together with the same hyperparameters and paired
initializations.
In this variant, the maximizer always retains the action
$\widetilde y_{t+1}$ proposed by the gradient update,
so $y_{t+1}=\widetilde y_{t+1}$ at every iteration.

\begin{figure}[h]
    \centering
    \includegraphics[width=0.8\linewidth]{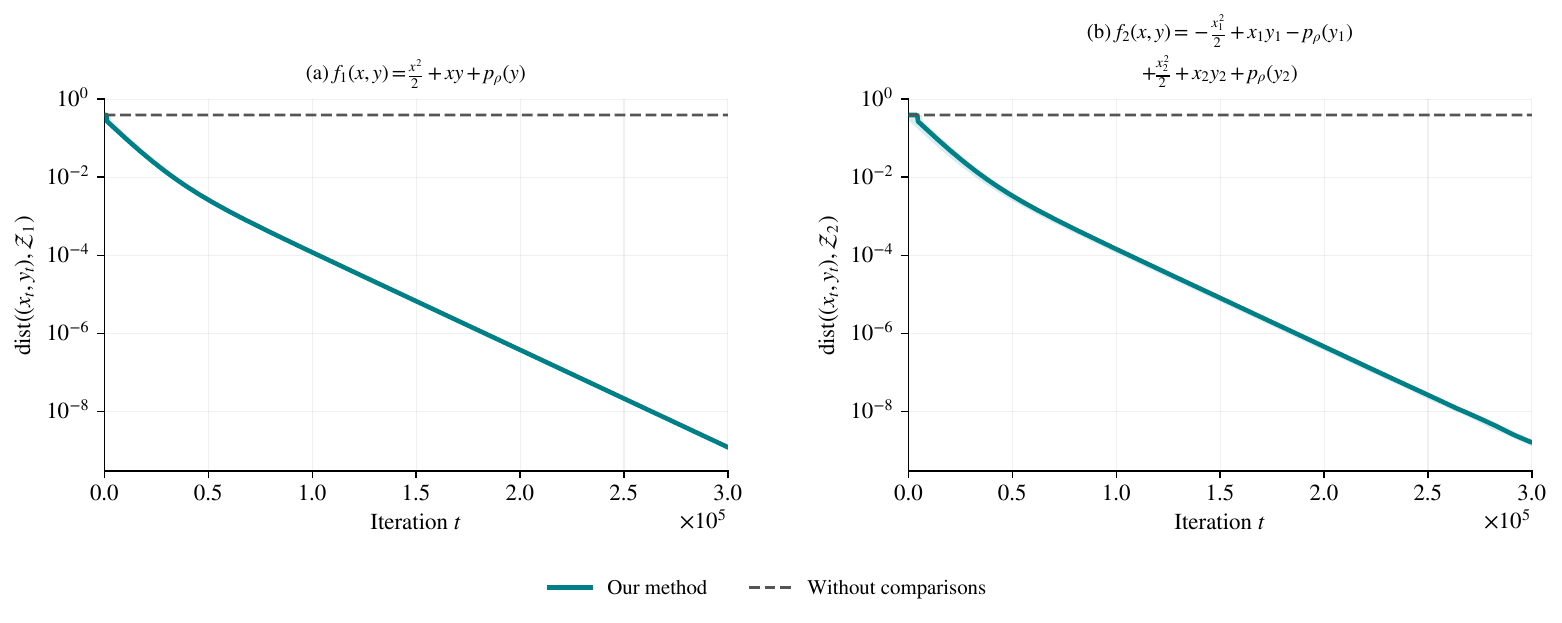}
    \caption{Effect of the comparison step on $f_1$
    (left) and $f_2$ (right).}
    \label{fig:comparison-ablation}
\end{figure}

Figure~\ref{fig:comparison-ablation} compares this variant
with our algorithm.
The left and right panels show distances to
$\mathcal Z_1$ for $f_1$ and $\mathcal Z_2$ for $f_2$,
respectively, over $300{,}000$ iterations.
Without comparisons, the iterates approach the incorrect
stationary point at the origin, leaving a positive
distance to the corresponding solution set.
With comparisons, our algorithm leaves the neighborhood
of the origin and approaches a calm local minimax solution.

In these two problems, the comparison step enables our
algorithm to leave the incorrect stationary attractor and
approach a calm local minimax solution. The same gradient
updates approach the incorrect point without comparisons.

\paragraph{Role of the comparison margin.}
The comparison step should detect improvements near an
incorrect stationary point while preserving convergence
near a valid solution.
To examine the role of the comparison margin, we consider
\begin{equation}
    f_3(x,y)=xy,
    \  x,y\in\RR,
    \label{eq:numerical-third-function}
\end{equation}
whose unique stationary point $(0,0)$ is a nondegenerate
calm local minimax point with a singular maximizer Hessian.

For this scalar objective, the random comparison step
proceeds as follows.
After the gradient update proposes the maximizer action
$\widetilde y_{t+1}$, our algorithm draws
$U_t$ uniformly from $[-1,1]$ and forms the candidate
\[
    \widehat y_{t+1}
    =\widetilde y_{t+1}+\ell U_t.
\]
It adopts this candidate only when
\[
    f_3(x_{t+1},\widehat y_{t+1})
    >
    f_3(x_{t+1},\widetilde y_{t+1})
    +\vartheta_t;
\]
otherwise, it retains $\widetilde y_{t+1}$.

We compare two versions of our algorithm.
The first one uses the adaptive margin design rule $\vartheta_t
    =
    \ell\sqrt{
        L\ell\left(
            \norm{\nabla f(x_t,y_t)}
            +s\norm{y_t-v_t}
        \right)
    }$
in~\eqref{eq:adaptive-margin}, and the other sets
$\vartheta_t=0$.
Both use $\eta=1/4$, $s=1/2$, $\tau=4$, $\ell=1/2$,
$v_0=y_0$, and the same 30 initializations satisfying $\norm{(x_0,y_0)}=10^{-2}\sqrt{2}$.
The two versions also use identical sequences of $U_t$, so only the comparison margin changes.

\begin{figure}[H]
    \centering
    \includegraphics[width=0.82\linewidth]{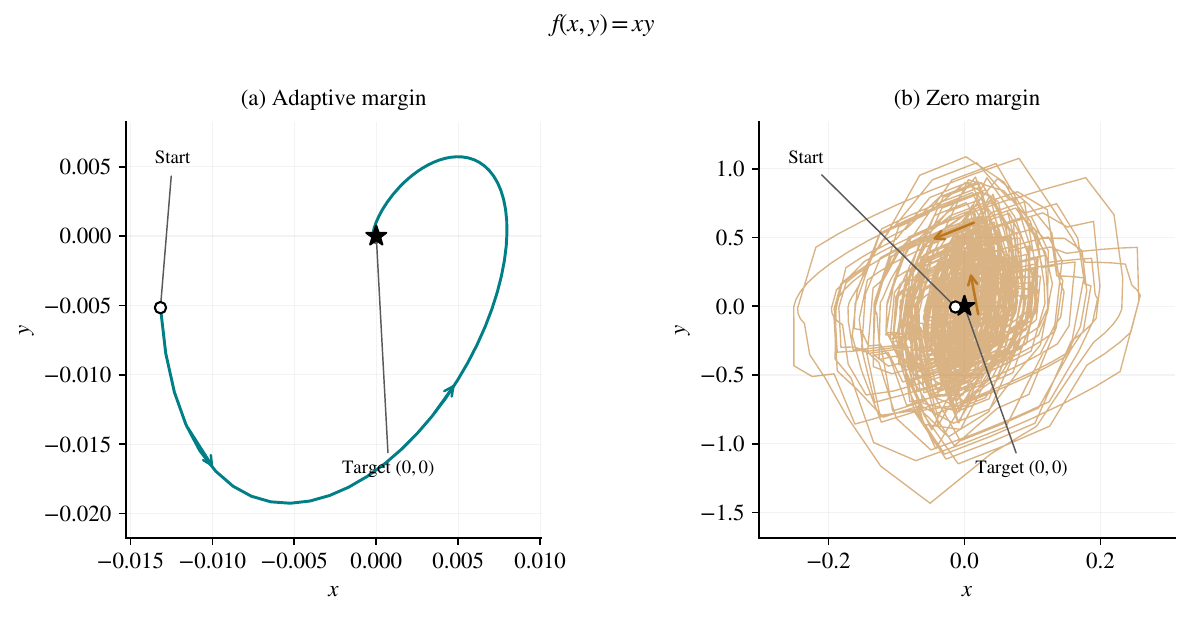}
    \caption{Trajectories for $f_3(x,y)=xy$:
    adaptive margin (left) and zero margin (right).}
    \label{fig:margin-phase}
\end{figure}

Figure~\ref{fig:margin-phase} displays representative
trajectories in the $(x,y)$-plane.
The left panel shows the adaptive-margin version, and the
right panel shows the zero-margin version.
The star marks the calm local minimax point $(0,0)$,
and the panels use different coordinate ranges to show
the complete trajectories.
With the adaptive margin, our algorithm approaches
the origin.
With zero margin, it continues making large moves
around the solution.

\begin{figure}[H]
    \centering
    \includegraphics[width=0.5\linewidth]{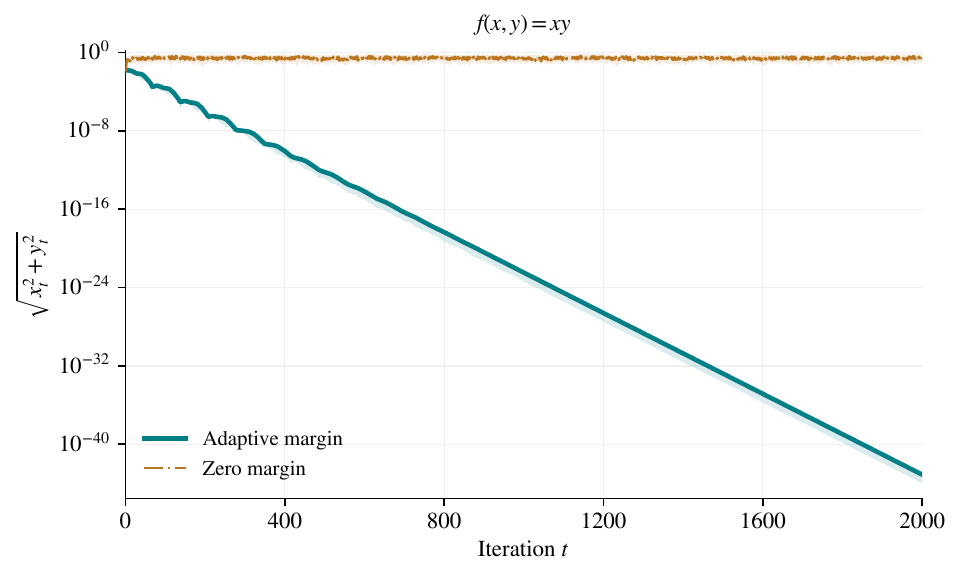}
    \caption{Distance to the origin for $f_3(x,y)=xy$
    with adaptive and zero margins.}
    \label{fig:margin-convergence}
\end{figure}

Figure~\ref{fig:margin-convergence} plots
$\sqrt{x_t^2+y_t^2}$ over $2{,}000$ iterations.
The adaptive-margin version consistently reduces
its distance to the origin, while the zero-margin
version remains away from it over the tested iterations.
Across the 30 runs, the adaptive-margin version accepts
no sampled candidates, whereas the zero-margin version
accepts a median of $1{,}005.5$ candidates.

This difference follows from the objective improvement
tested by the comparison oracle:
\[
\begin{aligned}
    &f_3(x_{t+1},\widehat y_{t+1})
    -f_3(x_{t+1},\widetilde y_{t+1})=x_{t+1}\ell U_t.
\end{aligned}
\]
Whenever $x_{t+1}\ne0$, a zero margin accepts the
candidate with probability $1/2$.
This remains true even when $x_{t+1}$ is arbitrarily
close to zero and the objective improvement is arbitrarily small.
However, the accepted change in $y$ is still $\ell U_t$,
whose magnitude need not decrease near the solution.
The adaptive margin rejects these small-gain moves
in the tested runs.
In these runs, the positive margin prevents the fixed-radius
comparison step from disrupting convergence near a valid
calm local minimax point.

\section{Discussion}
\label{sec:discussion}

In this paper, we develop an algorithm that combines
first-order gradient updates with zeroth-order function
comparisons to find nondegenerate calm local minimax points.
We establish local convergence to every such point from
sufficiently close initializations, for every allowed
sequence of random perturbations.
When every stationary point is nondegenerate and at least
one calm local minimax point exists, we further establish
almost-sure convergence of the entire iterate sequence
from any nonstationary initialization.
The global guarantee uses Gaussian random restarts to reach
a neighborhood where the local convergence result applies.
These results provide a partial resolution of the open
problem posed by \citet{chae2023open}.
The guarantees identify local minimax optima themselves
and include non-strict points with singular maximizer Hessians.

Several questions remain for future work.
The first is to develop first-order algorithms for general
local minimax points, including degenerate points and points
without the calmness property, toward a complete resolution
of the open problem.
A second direction is to establish finite-time oracle
complexity bounds that account for initialization and
restarts. Such bounds would quantify how dimension and local
geometry affect the total number of gradient and comparison
queries needed to attain an appropriate notion of approximate
local minimax optimality.
Finally, extending the method to stochastic gradients
and noisy function comparisons raises the question of
how to choose comparison margins and sample sizes that
account for estimation uncertainty while preserving
convergence guarantees.
\bibliographystyle{plainnat}
\bibliography{refs}

@article{lee2019first,
  title={First-order methods almost always avoid strict saddle points},
  author={Lee, Jason D. and Panageas, Ioannis and Piliouras, Georgios and Simchowitz, Max and Jordan, Michael I. and Recht, Benjamin},
  journal={Mathematical Programming},
  volume={176},
  pages={311--337},
  year={2019}
}

@inproceedings{jin2020local,
  title={What is local optimality in nonconvex-nonconcave minimax optimization?},
  author={Jin, Chi and Netrapalli, Praneeth and Jordan, Michael I.},
  booktitle={International Conference on Machine Learning},
  pages={4880--4889},
  year={2020}
}

@inproceedings{chae2023open,
  title={Open problem: Is there a first-order method that only converges to local minimax optima?},
  author={Chae, Jiseok and Kim, Kyuwon and Kim, Donghwan},
  booktitle={Conference on Learning Theory},
  pages={5957--5964},
  year={2023}
}

@inproceedings{chae2024two,
  title={Two-timescale extragradient for finding local minimax points},
  author={Chae, Jiseok and Kim, Kyuwon and Kim, Donghwan},
  booktitle={International Conference on Learning Representations},
  year={2024}
}

@inproceedings{kim2024double,
  title={Double-step alternating extragradient with increasing timescale separation for finding local minimax points: Provable improvements},
  author={Kim, Kyuwon and Kim, Donghwan},
  booktitle={International Conference on Machine Learning},
  pages={24059--24093},
  year={2024}
}

@article{ma2025calm,
  title={Calm local optimality for nonconvex-nonconcave minimax problems},
  author={Ma, Xiaoxiao and Yao, Wei and Ye, Jane J. and Zhang, Jin},
  journal={Set-Valued and Variational Analysis},
  volume={33},
  number={2},
  pages={13},
  year={2025}
}

@book{boyd1994linear,
  author={Boyd, Stephen and El Ghaoui, Laurent and Feron, Eric and Balakrishnan, Venkataramanan},
  title={Linear Matrix Inequalities in System and Control Theory},
  year={1994},
  series={Studies in Applied Mathematics},
  volume={15},
  publisher={SIAM},
}

@book{dontchev2009implicit,
  author={Dontchev, Asen L. and Rockafellar, R. Tyrrell},
  title={Implicit Functions and Solution Mappings: A View from Variational Analysis},
  series={Springer Monographs in Mathematics},
  publisher={Springer},
  year={2009},
}

@book{kellett2023introduction,
  author={Kellett, Christopher M. and Braun, Philipp},
  title={Introduction to Nonlinear Control: Stability, Control Design, and Estimation},
  publisher={Princeton University Press},
  year={2023},
  isbn={9780691240480}
}

@inproceedings{lee2016gradient,
  author={Lee, Jason D. and Simchowitz, Max and Jordan, Michael I. and Recht, Benjamin},
  title={Gradient Descent Only Converges to Minimizers},
  booktitle={Proceedings of the 29th Annual Conference on Learning Theory},
  series={Proceedings of Machine Learning Research},
  volume={49},
  pages={1246--1257},
  year={2016},
  publisher={PMLR},
}

@book{nemirovski1983problem,
  author={Nemirovsky, Arkadii S. and Yudin, David B.},
  title={Problem Complexity and Method Efficiency in Optimization},
  publisher={John Wiley \& Sons},
  year={1983},
  series={Wiley-Interscience Series in Discrete Mathematics},
  isbn={9780471103455}
}

@article{carmon2020lower,
  author={Carmon, Yair and Duchi, John C. and Hinder, Oliver and Sidford, Aaron},
  title={Lower Bounds for Finding Stationary Points {I}},
  journal={Mathematical Programming},
  volume={184},
  number={1--2},
  pages={71--120},
  year={2020},
}

@article{scheinberg2026functionfree,
  title={Function-free optimization via comparison oracles},
  author={Scheinberg, Katya and Xiong, Zikai},
  journal={arXiv preprint arXiv:2604.26867},
  year={2026}
}

@inproceedings{panageas2019first,
  title={First-order methods almost always avoid saddle points: The case of vanishing step-sizes},
  author={Panageas, Ioannis and Piliouras, Georgios and Wang, Xiao},
  booktitle={Advances in Neural Information Processing Systems},
  year={2019}
}

@inproceedings{fiez2021local,
  title={Local convergence analysis of gradient descent ascent with finite timescale separation},
  author={Fiez, Tanner and Ratliff, Lillian},
  booktitle={International Conference on Learning Representations},
  year={2021}
}

@inproceedings{fiez2020implicit,
  title={Implicit learning dynamics in {S}tackelberg games: Equilibria characterization, convergence analysis, and empirical study},
  author={Fiez, Tanner and Chasnov, Benjamin and Ratliff, Lillian},
  booktitle={International Conference on Machine Learning},
  year={2020}
}

@article{zhang2022optimality,
  title={Optimality and stability in non-convex smooth games},
  author={Zhang, Guojun and Poupart, Pascal and Yu, Yaoliang},
  journal={Journal of Machine Learning Research},
  volume={23},
  number={35},
  pages={1--71},
  year={2022}
}

@book{shub1987global,
  title={Global Stability of Dynamical Systems},
  author={Shub, Michael},
  year={1987},
  publisher={Springer-Verlag},
  address={New York}
}

@inproceedings{wang2020solving,
  title={On solving minimax optimization locally: A follow-the-ridge approach},
  author={Wang, Yuanhao and Zhang, Guodong and Ba, Jimmy},
  booktitle={International Conference on Learning Representations},
  year={2020}
}

@inproceedings{zhang2020single,
  title={A single-loop smoothed gradient descent-ascent algorithm for nonconvex-concave min-max problems},
  author={Zhang, Jiawei and Xiao, Peijun and Sun, Ruoyu and Luo, Zhi-Quan},
  booktitle={Advances in Neural Information Processing Systems},
  year={2020}
}

@inproceedings{zheng2023universal,
  author    = {Zheng, Taoli and Zhu, Linglingzhi and So, Anthony Man-Cho and Blanchet, Jos{\'e} and Li, Jiajin},
  title     = {Universal gradient descent ascent method for nonconvex-nonconcave minimax optimization},
  booktitle = {Advances in Neural Information Processing Systems},
  volume    = {36},
  year      = {2023}
}

@misc{zheng2024dissipative,
  author = {Zheng, Tianqi and Loizou, Nicolas and You, Pengcheng and Mallada, Enrique},
  title  = {Dissipative gradient descent ascent method: A control theory inspired algorithm for min-max optimization},
  year   = {2024},
  note   = {arXiv:2403.09090}
}

@book{wald1950statistical,
  author    = {Wald, Abraham},
  title     = {Statistical Decision Functions},
  publisher = {John Wiley \& Sons},
  address   = {New York},
  year      = {1950}
}

@book{bental2009robust,
  author    = {Ben-Tal, Aharon and El Ghaoui, Laurent and Nemirovski, Arkadi},
  title     = {Robust Optimization},
  series    = {Princeton Series in Applied Mathematics},
  publisher = {Princeton University Press},
  year      = {2009},
  isbn      = {9780691143682}
}

@book{basar1995hinfinity,
  author    = {Ba{\c{s}}ar, Tamer and Bernhard, Pierre},
  title     = {{$H_\infty$}-Optimal Control and Related Minimax Design Problems: A Dynamic Game Approach},
  edition   = {2nd},
  publisher = {Birkh{\"a}user},
  address   = {Boston, MA},
  year      = {1995}
}

@inproceedings{daskalakis2023stayon,
  title={{STay-ON-the-Ridge}: Guaranteed Convergence to
         Local Minimax Equilibrium in Nonconvex-Nonconcave Games},
  author={Daskalakis, Constantinos and Golowich, Noah and
          Skoulakis, Stratis and Zampetakis, Emmanouil},
  booktitle={Conference on Learning Theory},
  volume={195},
  series={Proceedings of Machine Learning Research},
  pages={5146--5198},
  year={2023},
  publisher={PMLR},
}

@book{nesterov2004introductory,
  title={Introductory Lectures on Convex Optimization:
         A Basic Course},
  author={Nesterov, Yurii},
  series={Applied Optimization},
  volume={87},
  publisher={Springer},
  year={2004},
}

@article{ghadimi2013stochastic,
  title={Stochastic First- and Zeroth-Order Methods
         for Nonconvex Stochastic Programming},
  author={Ghadimi, Saeed and Lan, Guanghui},
  journal={SIAM Journal on Optimization},
  volume={23},
  number={4},
  pages={2341--2368},
  year={2013},
}

@inproceedings{jin2017escape,
  title={How to Escape Saddle Points Efficiently},
  author={Jin, Chi and Ge, Rong and Netrapalli, Praneeth
          and Kakade, Sham M. and Jordan, Michael I.},
  booktitle={International Conference on Machine Learning},
  volume={70},
  series={Proceedings of Machine Learning Research},
  pages={1724--1732},
  year={2017},
  publisher={PMLR},
}

@inproceedings{hsieh2021limits,
  title={The Limits of Min-Max Optimization Algorithms:
         Convergence to Spurious Non-Critical Sets},
  author={Hsieh, Ya-Ping and Mertikopoulos, Panayotis
          and Cevher, Volkan},
  booktitle={International Conference on Machine Learning},
  volume={139},
  series={Proceedings of Machine Learning Research},
  pages={4337--4348},
  year={2021},
  publisher={PMLR},
}

@article{beck2003mirror,
  title={Mirror Descent and Nonlinear Projected Subgradient Methods for Convex Optimization},
  author={Beck, Amir and Teboulle, Marc},
  journal={Operations Research Letters},
  volume={31},
  number={3},
  pages={167--175},
  year={2003},
}

@article{parikh2014proximal,
  title={Proximal Algorithms},
  author={Parikh, Neal and Boyd, Stephen},
  journal={Foundations and Trends in Optimization},
  volume={1},
  number={3},
  pages={127--239},
  year={2014},
}

@inproceedings{reddi2016stochastic,
  title={Stochastic Variance Reduction for Nonconvex Optimization},
  author={Reddi, Sashank J. and Hefny, Ahmed and Sra, Suvrit and P{\'o}cz{\'o}s, Barnabas and Smola, Alex},
  booktitle={International Conference on Machine Learning},
  volume={48},
  series={Proceedings of Machine Learning Research},
  pages={314--323},
  year={2016},
  publisher={PMLR},
}

@inproceedings{fang2018spider,
  title={{SPIDER}: Near-Optimal Non-Convex Optimization via Stochastic Path-Integrated Differential Estimator},
  author={Fang, Cong and Li, Chris Junchi and Lin, Zhouchen and Zhang, Tong},
  booktitle={Advances in Neural Information Processing Systems},
  volume={31},
  year={2018},
}

@article{nemirovski2004prox,
  title={Prox-Method with Rate of Convergence {$O(1/t)$} for Variational Inequalities with Lipschitz Continuous Monotone Operators and Smooth Convex-Concave Saddle Point Problems},
  author={Nemirovski, Arkadi},
  journal={SIAM Journal on Optimization},
  volume={15},
  number={1},
  pages={229--251},
  year={2004},
}

@inproceedings{daskalakis2018training,
  title={Training {GANs} with Optimism},
  author={Daskalakis, Constantinos and Ilyas, Andrew and Syrgkanis, Vasilis and Zeng, Haoyang},
  booktitle={International Conference on Learning Representations},
  year={2018},
}

@inproceedings{mokhtari2020unified,
  title={A Unified Analysis of Extra-gradient and Optimistic Gradient Methods for Saddle Point Problems: Proximal Point Approach},
  author={Mokhtari, Aryan and Ozdaglar, Asuman and Pattathil, Sarath},
  booktitle={International Conference on Artificial Intelligence and Statistics},
  volume={108},
  series={Proceedings of Machine Learning Research},
  pages={1497--1507},
  year={2020},
  publisher={PMLR},
}

@article{nesterov2017random,
  title={Random Gradient-Free Minimization of Convex Functions},
  author={Nesterov, Yurii and Spokoiny, Vladimir},
  journal={Foundations of Computational Mathematics},
  volume={17},
  number={2},
  pages={527--566},
  year={2017},
}

@article{duchi2015optimal,
  title={Optimal Rates for Zero-Order Convex Optimization: The Power of Two Function Evaluations},
  author={Duchi, John C. and Jordan, Michael I. and Wainwright, Martin J. and Wibisono, Andre},
  journal={IEEE Transactions on Information Theory},
  volume={61},
  number={5},
  pages={2788--2806},
  year={2015},
}

@article{bergou2020stochastic,
  title={Stochastic Three Points Method for Unconstrained Smooth Minimization},
  author={Bergou, El Houcine and Gorbunov, Eduard and Richt{\'a}rik, Peter},
  journal={SIAM Journal on Optimization},
  volume={30},
  number={4},
  pages={2726--2749},
  year={2020},
}

@inproceedings{jamieson2012query,
  title={Query Complexity of Derivative-Free Optimization},
  author={Jamieson, Kevin G. and Nowak, Robert D. and Recht, Benjamin},
  booktitle={Advances in Neural Information Processing Systems},
  volume={25},
  pages={2672--2680},
  year={2012},
}

@inproceedings{karabag2021smooth,
  title={Smooth Convex Optimization Using Sub-Zeroth-Order Oracles},
  author={Karabag, Mustafa O. and Neary, Cyrus and Topcu, Ufuk},
  booktitle={AAAI Conference on Artificial Intelligence},
  volume={35},
  pages={3815--3822},
  year={2021},
}

@inproceedings{saha2021dueling,
  title={Dueling Convex Optimization},
  author={Saha, Aadirupa and Koren, Tomer and Mansour, Yishay},
  booktitle={International Conference on Machine Learning},
  volume={139},
  series={Proceedings of Machine Learning Research},
  pages={9245--9254},
  year={2021},
  publisher={PMLR},
}

@inproceedings{saha2025dueling,
  title={Dueling Convex Optimization with General Preferences},
  author={Saha, Aadirupa and Koren, Tomer and Mansour, Yishay},
  booktitle={International Conference on Machine Learning},
  volume={267},
  series={Proceedings of Machine Learning Research},
  pages={52552--52564},
  year={2025},
  publisher={PMLR},
}

@article{cai2022onebit,
  title={A One-Bit, Comparison-Based Gradient Estimator},
  author={Cai, HanQin and McKenzie, Daniel and Yin, Wotao and Zhang, Zhenliang},
  journal={Applied and Computational Harmonic Analysis},
  volume={60},
  pages={242--266},
  year={2022},
}

@article{hu2026optimizing,
  title={Optimizing the Preconditioner: A Black-box Online-to-Nonconvex Conversion with Static Regret Minimization Oracles},
  author={Hu, Haichen and Simchi-Levi, David},
  journal={arXiv preprint arXiv:2607.17607},
  year={2026}
}

@article{ahn2024general,
  title={General framework for online-to-nonconvex conversion: Schedule-free SGD is also effective for nonconvex optimization},
  author={Ahn, Kwangjun and Magakyan, Gagik and Cutkosky, Ashok},
  journal={arXiv preprint arXiv:2411.07061},
  year={2024}
}

@inproceedings{cutkosky2023optimal,
  title={Optimal stochastic non-smooth non-convex optimization through online-to-non-convex conversion},
  author={Cutkosky, Ashok and Mehta, Harsh and Orabona, Francesco},
  booktitle={International Conference on Machine Learning},
  pages={6643--6670},
  year={2023},
  organization={PMLR}
}

\appendix
\section{Mathematical Tools}
\label{sec:mathematical-tools}

In this section, we collect the mathematical tools used
to establish the convergence guarantees in
Section~\ref{sec:theory}.
These results concern quadratic forms on subspaces,
local error bounds, and the behavior of iterations
near stable and unstable fixed points.

We first recall a form of Finsler's lemma
\citep{boyd1994linear}.
It relates positive definiteness on a subspace to
positive definiteness after adding a quadratic penalty.

\begin{lemma}[Finsler's lemma]
\label{lem:tools-finsler}
Let $Q\in\RR^{m\times m}$ be symmetric and let
$D\in\RR^{m\times r}$.
The following conditions are equivalent:
\begin{enumerate}[label=(\alph*)]
    \item
    $u^\top Q u>0$ for every nonzero
    $u\in\ker D^\top$.
    \item
    There exists $\lambda>0$ such that
    \[
        Q+\lambda DD^\top\succ0.
    \]
\end{enumerate}
No rank assumption on $D$ is required.
If $r=0$, we interpret $DD^\top=0$.
\end{lemma}

The next result is a direct consequence of the inverse
function theorem
\citep[Theorem~1A.1]{dontchev2009implicit}.
It bounds the distance to a nondegenerate zero
in terms of the residual.

\begin{lemma}[\cite{dontchev2009implicit}]
\label{lem:tools-error-bound}
Let $\mathcal U\subseteq\RR^m$ be open, and let
$G:\mathcal U\to\RR^m$ be continuously differentiable.
Suppose that $u^*\in\mathcal U$, $G(u^*)=0$, and
$DG(u^*)$ is nonsingular.
Then there exist constants $r>0$ and $0<c\le C$ such that
$\{u:\norm{u-u^*}<r\}\subseteq\mathcal U$ and
\begin{equation}
    c\norm{u-u^*}
    \le
    \norm{G(u)}
    \le
    C\norm{u-u^*}
    \ \text{whenever }\norm{u-u^*}<r.
    \label{eq:tools-error-bound}
\end{equation}
In particular, $u^*$ is the unique zero of $G$
in this neighborhood.
\end{lemma}

We next recall the discrete-time linearization criterion
for local exponential stability
\citep[Theorem~5.9]{kellett2023introduction}.
For a square matrix $A$, let $\rho(A)$ denote its
spectral radius, namely the largest modulus of its
eigenvalues.

\begin{lemma}[\cite{kellett2023introduction}]
\label{lem:tools-linear-stability}
Let $T:\RR^m\to\RR^m$ be continuously differentiable.
Suppose that $T(u^*)=u^*$ and
\[
    \rho(DT(u^*))<1.
\]
Then there exist constants $r>0$, $C\ge1$, and
$q\in(0,1)$ such that
\begin{equation}
    \norm{T^n(u)-u^*}
    \le
    Cq^n\norm{u-u^*},
    \  n=0,1,2,\ldots,
    \label{eq:tools-linear-stability}
\end{equation}
for every $u$ satisfying $\norm{u-u^*}<r$.
Here $T^n$ denotes the $n$-fold composition of $T$,
with $T^0$ equal to the identity map.
\end{lemma}

Finally, we recall the center-stable manifold theorem
\citep[Theorem~III.7]{shub1987global}.
It describes the set of initial points whose iterates
remain near a fixed point. We use it to exclude
convergence to unstable stationary points.

\begin{lemma}[\cite{shub1987global}]
\label{lem:tools-center-stable}
Let $\mathcal U\subseteq\RR^m$ be open, and let
$T:\mathcal U\to\RR^m$ be continuously differentiable,
with $DT(u)$ nonsingular for every $u\in\mathcal U$.
Suppose that $T(u^*)=u^*$ for some $u^*\in\mathcal U$.

Let $E^{cs}$ denote the real invariant subspace
associated with the generalized eigenspaces of
$DT(u^*)$ whose eigenvalues have modulus at most one.
Then there exist a neighborhood
$V\subseteq\mathcal U$ of $u^*$ and a $C^1$
embedded submanifold $W_{\mathrm{loc}}^{cs}$
containing $u^*$ such that
\[
    T_{u^*}W_{\mathrm{loc}}^{cs}=E^{cs}
\]
and
\begin{equation}
    \bigl\{
        u\in V:
        T^n(u)\in V
        \text{ for every }n\ge0
    \bigr\}
    \subseteq W_{\mathrm{loc}}^{cs}.
    \label{eq:tools-center-stable}
\end{equation}
Here $T_{u^*}W_{\mathrm{loc}}^{cs}$ denotes the tangent
space of the manifold at $u^*$.

If $DT(u^*)$ has an eigenvalue with modulus greater
than one, then
$\dim W_{\mathrm{loc}}^{cs}<m$.
\end{lemma}

\section{Proofs in Section \ref{sec:theory}}
\label{app:proofs_theory}

In this section, we prove the local and global convergence
guarantees for Algorithm~\ref{alg:hybrid}.
We first establish local convergence of the gradient and
comparison updates with fixed parameters. We then show how
the restart procedure selects suitable parameters without
knowing the target point. For global convergence, we also
prove that the algorithm almost surely avoids stationary
limits that are not calm local minimax points and makes
only finitely many restarts.
We combine these lemmas to prove
Theorems~\ref{thm:local-convergence}
and~\ref{thm:global-convergence}.
The proofs of the lemmas are collected in
Appendix~\ref{app:proofs_app_theory}.
Throughout this section, Assumption~\ref{ass:smooth} holds,
and the algorithm is run without a finite iteration budget,
as specified in Section~\ref{sec:theory}.

We first express calm local minimax optimality through
the curvature condition used to establish local convergence.
At a point $z^*=(x^*,y^*)$, write
\begin{equation}
    H=\frac1L\nabla^2f(z^*)
      =\begin{pmatrix}A&C\\C^\top&B\end{pmatrix}.
    \label{eq:hessian-blocks}
\end{equation}
Whenever $B\preceq0$, choose orthonormal $y$ coordinates and
an orthonormal basis $P$ of $\ker C_2^\top$ such that
\begin{equation}
\begin{aligned}
    B&=\begin{pmatrix}-D&0\\0&0\end{pmatrix},
    &C&=\begin{pmatrix}C_1&C_2\end{pmatrix},\  D\succ0,\\
    S&=A+C_1D^{-1}C_1^\top,
    &S_0&=P^\top SP.
\end{aligned}
\label{eq:reduced-curvature}
\end{equation}
Empty blocks are omitted; their determinants are one.

\begin{lemma}
\label{lem:curvature-characterization}
For the matrices in
\eqref{eq:hessian-blocks}--\eqref{eq:reduced-curvature},
assume $\det H\ne0$ and $B\preceq0$. Then
\[
    \operatorname{rank}C_2=\dim\ker B,
    \  \det S_0\ne0.
\]
If $\nabla f(z^*)=0$ and $y^*$ locally
maximizes $f(x^*,\cdot)$, then
\[
    z^*\text{ is a calm local minimax point}
    \ \Longleftrightarrow\  S_0\succ0.
\]
\end{lemma}

The condition $S_0\succ0$ will ensure stability of the
gradient update. The next lemma establishes the corresponding
eigenvalue bound for its linearization.

\begin{lemma}
\label{lem:spectral-stability}
Assume $\norm{H}\le1$, $\det H\ne0$, $B\preceq0$,
and $S_0\succ0$. For $\gamma>0$, define
\begin{equation}
    \mathsf X_\gamma=
    \begin{pmatrix}
        \gamma^2A&\gamma C&0\\
        -\gamma C^\top&\gamma I-B&-\gamma^2I\\
        0&-\gamma^2I&\gamma^3I
    \end{pmatrix}.
    \label{eq:spectral-matrix}
\end{equation}
There exists $\gamma_*\in(0,1)$ such that
\[
    0<\gamma<\gamma_*
    \ \longrightarrow\ 
    |1-\lambda/4|<1
    \ \text{for every eigenvalue }\lambda
    \text{ of }\mathsf X_\gamma.
\]
\end{lemma}

We next combine this eigenvalue bound with the comparison rule.
Near the target, the adaptive margin makes every comparison
retain the gradient proposal, so the local convergence bound
also holds for the combined updates.

\begin{lemma}
\label{lem:fixed-parameter-convergence}
Suppose Assumption~\ref{ass:smooth} holds.
Let $z^*=(x^*,y^*)$ be a calm local minimax point with
$\nabla^2 f(z^*)$ nonsingular, and write
$w^*=(x^*,y^*,y^*)$.
There exists $\gamma_*\in(0,1)$ such that, for every
$\gamma\in(0,\gamma_*)$, the parameters
\[
    s=L\gamma,\  \tau=\gamma^{-2},\ 
    \ell=\gamma,\  \eta=\frac{1}{4L}
\]
have the following property.
There exist $r>0$, $C\ge1$, and $q\in(0,1)$ such that,
for every $w_0$ satisfying $\norm{w_0-w^*}<r$, the updates
\eqref{eq:base-update} and \eqref{eq:comparison-update},
applied without the restart check, satisfy
\begin{equation}
    \norm{w_t-w^*}\le Cq^t\norm{w_0-w^*},
    \  t=0,1,2,\ldots.
    \label{eq:fixed-parameter-rate}
\end{equation}
Every comparison query returns zero.
Both conclusions hold for every sequence $(U_t)_{t\ge0}$
satisfying $\norm{U_t}\le1$.
\end{lemma}

The preceding result requires suitable fixed parameters.
We next show that the restart check guarantees convergence
to a stationary point whenever a trial continues forever,
even without knowing whether its parameters are suitable.

\begin{lemma}
\label{lem:continuing-trial}
Suppose Assumption~\ref{ass:smooth} holds.
Consider a trial $j$ with $R_j>0$ for which
\eqref{eq:restart-check} holds at every step.
There exists $z^\infty=(x^\infty,y^\infty)$ such that
\begin{equation}
\begin{aligned}
    \nabla f(z^\infty)&=0,\\
    \norm{w_{t_j+n}-(x^\infty,y^\infty,y^\infty)}
    &\le\alpha_j(1-\alpha_j)^n,
    \  n=0,1,2,\ldots.
\end{aligned}
\label{eq:continuing-trial-rate}
\end{equation}
The conclusion holds for every sequence $(U_t)$ satisfying
$\norm{U_t}\le1$.
\end{lemma}

We must also ensure that some trial satisfies every restart check.
The next lemma shows that this holds near the target point
by combining the local convergence bound, initialization
perturbation, and decreasing step bound.

\begin{lemma}
\label{lem:successful-trial}
Suppose Assumption~\ref{ass:smooth} holds.
Let $z^*=(x^*,y^*)$ be calm local minimax with
$\det\nabla^2 f(z^*)\ne0$, and set $w^*=(x^*,y^*,y^*)$.
There exists $\gamma_*\in(0,1)$ such that, for each
$\gamma\in(0,\gamma_*)$, there is $\delta>0$ with the
following property.
For any trial $j$ of Algorithm~\ref{alg:hybrid} satisfying
\[
    \gamma_j=\gamma,\ \norm{\bar w_j-w^*}<\delta,
\]
if $R_j=0$, the algorithm returns $z^*$.
If $R_j>0$, then
\begin{equation}
\begin{aligned}
    \mathcal C_f\bigl(
        x_{t+1}^{\mathrm{cand}};
        \widetilde y_{t+1},\widehat y_{t+1},\vartheta_t
    \bigr)&=0,\  \norm{w_{t+1}^{\mathrm{cand}}-w_t}
    &\le\alpha_j^2(1-\alpha_j)^{t-t_j},
    \ t\ge t_j,\\
    w_{t_j+n}&\longrightarrow w^*,
    \ n\to\infty.
\end{aligned}
\label{eq:successful-trial-conclusions}
\end{equation}
These conclusions hold for every $Z_j$ and every sequence
$(U_t)$ satisfying $\norm{Z_j}\le1$ and $\norm{U_t}\le1$.
\end{lemma}

Before reaching a suitable scale, earlier trials may fail.
The next lemma ensures that any fixed number of these trials
remain near the target when the initial point is sufficiently close.

\begin{lemma}
\label{lem:local-restarts}
Suppose Assumption~\ref{ass:smooth} holds.
Let $\nabla f(z^*)=0$ and $w^*=(x^*,y^*,y^*)$.
For every integer $J\ge1$ and $\rho>0$, there exists $r>0$
such that the initialization conditions
\[
    w_0=(x_0,y_0,y_0),\ 
    \norm{w_0-w^*}<r,\ 
    \nabla f(x_0,y_0)\ne0
\]
imply, for every reached trial $j\le\min\{J,K\}$,
\begin{equation}
\begin{aligned}
    \norm{\bar w_j-w^*}&<\rho,\\
    \norm{w_t-w^*}&<\rho
    \ \text{for every state $w_t$ produced in trial }j.
\end{aligned}
\label{eq:finite-trial-neighborhood}
\end{equation}
The conclusions hold for every choice of $(Z_j)_{j\ge1}$
and $(U_t)_{t\ge0}$ satisfying
$\norm{Z_j}\le1$ and $\norm{U_t}\le1$.
\end{lemma}

These lemmas suffice for the local convergence theorem.
For global convergence, we must also exclude incorrect stationary
limits; we first show that the randomized updates preserve
absolute continuity.

\begin{lemma}
\label{lem:absolute-continuity}
Suppose Assumption~\ref{ass:smooth} holds.
Fix $s\in(0,L]$, $\tau\ge1$, $\ell>0$, and $\eta=1/(4L)$.
Apply \eqref{eq:base-update}--\eqref{eq:comparison-update}
without restarts, with independent random vectors $U_t$, each
uniform on the unit ball in $\RR^{d_2}$.
If $w_0$ has an absolutely continuous law and is independent
of $(U_t)_{t\ge0}$, then, for every Borel set
$N\subseteq\RR^{d_1+2d_2}$, we have
\[
    \operatorname{Leb}(N)=0
    \ \longrightarrow\ 
    \Pr(w_t\in N)=0,\  t=0,1,2,\ldots.
\]
\end{lemma}

Absolute continuity allows a measure-zero argument to exclude
incorrect limits. We next identify an unstable
direction when local maximality in $y$ holds but calm local
minimax optimality fails.

\begin{lemma}
\label{lem:gradient-instability}
Let $z^*$ be a nondegenerate stationary point that is not calm
local minimax, and suppose $y^*$ locally maximizes $f(x^*,\cdot)$.
For every $s\in(0,L]$, $\tau\ge1$, and $\eta=1/(4L)$,
the Jacobian matrix $\mathsf J$ of the three updates in
\eqref{eq:base-update} at $w^*=(x^*,y^*,y^*)$ satisfies
\[
    \mathsf J a=\lambda a
    \ \text{for some }a\ne0\text{ and }\lambda>1.
\]
\end{lemma}

The unstable direction excludes convergence in the preceding
case through a measure-zero argument. When local maximality
in $y$ fails instead, the random comparisons detect an improvement;
the next lemma combines these two cases.

\begin{lemma}
\label{lem:avoidance-fixed-parameters}
Suppose Assumption~\ref{ass:smooth} holds.
Fix $s\in(0,L]$, $\tau\ge1$, and $\ell>0$, and set
$\eta=1/(4L)$.
Apply \eqref{eq:base-update} and
\eqref{eq:comparison-update} without the restart check.
Suppose $w_0$ has a distribution absolutely continuous with
respect to Lebesgue measure on $\RR^{d_1+2d_2}$ and is
independent of $(U_t)_{t\ge0}$, whose elements are independent
and uniform on the unit ball in $\RR^{d_2}$.
Then the distribution of $w_t$ is absolutely continuous
for every finite $t$.
Moreover, for every nondegenerate stationary point
$z^*=(x^*,y^*)$ that is not a calm local minimax point,
\begin{equation}
    \Pr\!\left(
        \lim_{t\to\infty}w_t=(x^*,y^*,y^*)
    \right)=0.
    \label{eq:avoid-incorrect-limit}
\end{equation}
\end{lemma}

To apply this avoidance result after a restart, we use the fresh
random perturbation of each trial center. Under global
nondegeneracy, stationary points form a countable set,
so we can exclude all incorrect stationary limits simultaneously.

\begin{lemma}
\label{lem:correct-trial-limits}
Suppose Assumption~\ref{ass:smooth} holds, every stationary
point is nondegenerate, and $\nabla f(x_0,y_0)\ne0$.
With probability one, Algorithm~\ref{alg:hybrid} satisfies
\[
    R_j>0\ \text{for every reached trial }j.
\]
On the same probability-one event, each trial $j$ that
continues forever has a limit satisfying
\begin{equation}
\begin{gathered}
    w_{t_j+n}\longrightarrow(x^\infty,y^\infty,y^\infty),\ (x^\infty,y^\infty)\text{ is a calm local minimax point}.
\end{gathered}
\label{eq:correct-final-trial}
\end{equation}
\end{lemma}

It remains to rule out infinitely many restarts.
The next lemma uses the positive probability of sampling
a center and scale covered by Lemma~\ref{lem:successful-trial}
to show that only finitely many trials fail.

\begin{lemma}
\label{lem:global-restarts}
Suppose Assumption~\ref{ass:smooth} holds and $f$ has a
nondegenerate calm local minimax point.
For every initialization with $\nabla f(x_0,y_0)\ne0$,
Algorithm~\ref{alg:hybrid} satisfies
\begin{equation}
    \Pr\!\left(
        \bigl|\{j:\text{trial }j\text{ fails}\}\bigr|<\infty
    \right)=1.
    \label{eq:finite-restarts}
\end{equation}
\end{lemma}

We now prove local convergence by keeping the first trials near
the target and applying Lemma~\ref{lem:successful-trial}
once a sufficiently small scale is reached.

\begin{proof}[Proof of Theorem \ref{thm:local-convergence}]
Fix any sequences $(Z_j)_{j\ge1}$ and $(U_t)_{t\ge0}$
satisfying $\norm{Z_j}\le1$ and $\norm{U_t}\le1$.
Choose an integer $J_*\ge1$ such that $2^{-J_*}$ is
smaller than the threshold in Lemma~\ref{lem:successful-trial},
and let $\delta>0$ be the corresponding radius for trial centers.
Since $\nabla^2 f(z^*)$ is nonsingular,
Lemma~\ref{lem:tools-error-bound}, applied to $\nabla f$,
gives $r_*>0$ such that
\begin{equation}
    \nabla f(z)=0,\  \norm{z-z^*}<r_*
    \ \longrightarrow\  z=z^*.
    \label{eq:local-isolated-stationary}
\end{equation}
Fix $0<\rho<\min\{\delta,r_*/2\}$, and let $r_0>0$
be supplied by Lemma~\ref{lem:local-restarts} with
$J=J_*$ and this $\rho$.
Choose
\[
    0<r<\min\left\{r_*,\frac{r_0}{\sqrt{2}},J_*^{-3}\right\}.
\]
Every allowed initialization satisfies $    \norm{w_0-w^*}^2
    =\norm{x_0-x^*}^2+2\norm{y_0-y^*}^2
    \le2\norm{z_0-z^*}^2<r_0^2.$

If $\nabla f(z_0)=0$, \eqref{eq:local-isolated-stationary}
shows that the initial return is $z^*$.
Otherwise, the Hessian bound and $\nabla f(z^*)=0$ give
\begin{align*}
    0<\alpha_{\mathrm{init}}
    &\le\left(\frac{\norm{\nabla f(z_0)}}{L}\right)^{1/3}
     \le\norm{z_0-z^*}^{1/3}<\frac{1}{J_*},\ K=\left\lceil\alpha_{\mathrm{init}}^{-1}\right\rceil
       \ge J_*.
\end{align*}
Lemma~\ref{lem:local-restarts} therefore places the centers
and accepted states of every reached trial $j\le J_*$
inside the $\rho$-ball around $w^*$.

If such a trial has $R_j=0$, then
\[
    \nabla f(\bar x_j,\bar y_j)=0,\ 
    \bar y_j=\bar v_j,\ 
    \norm{(\bar x_j,\bar y_j)-z^*}<\rho<r_*.
\]
Equation~\eqref{eq:local-isolated-stationary} identifies
the returned point as $z^*$.
If one of these trials continues forever,
Lemma~\ref{lem:continuing-trial} gives a stationary limit
$z^\infty$ satisfying $    \norm{z^\infty-z^*}\le\rho<r_*$. Again \eqref{eq:local-isolated-stationary} yields
$z^\infty=z^*$.

It remains to consider the case in which the first
$J_*-1$ trials all fail.
Trial $J_*$ then has scale $2^{-J_*}$ and center satisfying
\[
    \norm{\bar w_{J_*}-w^*}<\rho<\delta.
\]
Lemma~\ref{lem:successful-trial} shows that this trial
either returns $z^*$ or continues forever and converges
to $w^*$.
There are consequently at most $J_*-1$ restarts, and every
trial uses \eqref{eq:local-trial-center} because $J_*\le K$.
For the final trial, \eqref{eq:continuing-trial-rate}
is exactly \eqref{eq:local-convergence-rate}.
\end{proof}

For global convergence, Lemma~\ref{lem:global-restarts}
gives finitely many restarts, and
Lemma~\ref{lem:correct-trial-limits} identifies the final limit.

\begin{proof}[Proof of Theorem \ref{thm:global-convergence}]
The assumed calm local minimax point is nondegenerate
because every stationary point is nondegenerate.
Lemma~\ref{lem:global-restarts} therefore implies that,
with probability one, the algorithm terminates or begins
a final trial after finitely many restarts.
The initial point is nonstationary, and
Lemma~\ref{lem:correct-trial-limits} gives
\[
    \Pr\!\left(R_j>0\text{ for every reached trial }j\right)=1.
\]
Thus, almost surely, no stationarity check terminates the run.
There is a final trial $J$, and the same lemma gives a
calm local minimax point $z^\infty=(x^\infty,y^\infty)$
such that
\[
    w_{t_J+n}\longrightarrow
    w^\infty=(x^\infty,y^\infty,y^\infty).
\]
Every preceding trial failed at a finite step, and there
are only finitely many such trials. Hence $t_J<\infty$,
and convergence of this tail implies convergence of the full
sequence $w_t$.
Finally, Lemma~\ref{lem:continuing-trial} yields
\[
    \norm{w_{t_J+n}-w^\infty}
    \le\alpha_J(1-\alpha_J)^n,
    \  n=0,1,2,\ldots,
\]
which proves \eqref{eq:global-convergence-rate}.
\end{proof}

\section{Proofs of Lemmas in Section \ref{app:proofs_theory}}
\label{app:proofs_app_theory}

In this section, we prove the lemmas in their stated order.

\begin{proof}[Proof of Lemma \ref{lem:curvature-characterization}]
In the coordinates of \eqref{eq:reduced-curvature},
\[
    H=\begin{pmatrix}
        A&C_1&C_2\\ C_1^\top&-D&0\\ C_2^\top&0&0
    \end{pmatrix}.
\]
If $C_2b=0$, then
\[
    H\begin{pmatrix}0\\0\\b\end{pmatrix}=0
    \ \longrightarrow\  b=0.
\]
Thus $\operatorname{rank}C_2=\dim\ker B$.
If $S_0a=0$, set $u=Pa$. Then
\[
    C_2^\top u=0,\  P^\top Su=0,
    \  Su\in\operatorname{range}C_2.
\]
Choose $b$ with $C_2b=-Su$. Direct multiplication gives
\[
    H\begin{pmatrix}u\\D^{-1}C_1^\top u\\b\end{pmatrix}
    =\begin{pmatrix}Su+C_2b\\0\\0\end{pmatrix}=0.
\]
Hence $u=0$ and $a=P^\top u=0$, proving $\det S_0\ne0$.

Now suppose that $z^*$ is stationary and $y^*$ locally
maximizes $f(x^*,\cdot)$.
Suppose first that $z^*$ is calm, with response constant $M$.
For a unit vector $u\in\ker C_2^\top$ and $h=(h_1,h_2)$
with $\norm{h}\le M\delta$, Taylor expansion gives
\begin{align*}
    \frac{f(x^*+\delta u,y^*+h)-f(z^*)}{L}
    &=\frac{\delta^2}{2}u^\top Au
       +\delta u^\top C_1h_1-\frac12h_1^\top Dh_1
       +o(\delta^2)\\
    &=\frac{\delta^2}{2}u^\top Su
       -\frac12\norm{D^{1/2}h_1
                    -\delta D^{-1/2}C_1^\top u}^2
       +o(\delta^2).
\end{align*}
The remainder is uniform over $\norm{h}\le M\delta$.
The calm inequality therefore implies
\[
    0\le\frac{\delta^2}{2}u^\top Su+o(\delta^2)
    \ \longrightarrow\  S_0\succeq0.
\]
Since $\det S_0\ne0$, we obtain $S_0\succ0$.

Conversely, suppose $S_0\succ0$.
Lemma~\ref{lem:tools-finsler} supplies $\lambda>0$ such that
\[
    S+\lambda C_2C_2^\top\succ0.
\]
For a minimizer displacement $a$, choose
\[
    h(a)=\begin{pmatrix}
        D^{-1}C_1^\top a\\(\lambda/2)C_2^\top a
    \end{pmatrix},
    \ \norm{h(a)}\le M\norm{a}
\]
with $M>0$ larger than the norm of this linear operator.
Then
\begin{align*}
    \frac{f(x^*+a,y^*+h(a))-f(z^*)}{L}
    &=\frac12a^\top(S+\lambda C_2C_2^\top)a
       +o(\norm{a}^2)\ge0
\end{align*}
for all sufficiently small $a$.
Choose $\rho>0$ within the neighborhood of local maximality.
For small $\delta\le\rho$ and $\norm{a}\le\delta$, this yields
\[
    f(x^*,y)\le f(z^*)
    \le f(x^*+a,y^*+h(a))
    \le\max_{\norm{y'-y^*}\le M\delta}f(x^*+a,y')
\]
whenever $\norm{y-y^*}\le\delta$, which is
\eqref{eq:calm-definition}.
\end{proof}

\begin{proof}[Proof of Lemma \ref{lem:spectral-stability}]
Put $k=\dim\ker B$, $r=d_2-k$, $m=d_1-k$, and $d=d_1+2d_2$,
and denote the eigenvalues of $\mathsf X_\gamma$ by
$\lambda_1(\gamma),\ldots,\lambda_d(\gamma)$.
Lemma~\ref{lem:curvature-characterization} gives
$\operatorname{rank}C_2=k$.
Complete $P$ to orthonormal $x$ coordinates $(u,\xi)$ so that
\[
    C_2=\begin{pmatrix}0\\E\end{pmatrix},
    \  E\in\RR^{k\times k},\  \det E\ne0.
\]
Write $y=(p,q)$ and $v=(a,b)$, with block sizes $(r,k)$.
Let $\delta_1,\ldots,\delta_r$ denote the eigenvalues of the block $D$
within this proof. The only nonzero block of $\mathsf X_0$ is $D$, so continuity
of polynomial roots gives $\lambda_i(\gamma)\to\delta_i$
for $1\le i\le r$ and $\lambda_i(\gamma)\to0$ otherwise.

For complex $\lambda$ near zero, put
\[
    R=(D+\gamma I-\lambda I)^{-1},\ 
    S_{\gamma,\lambda}=A+C_1RC_1^\top.
\]
Eliminating the $p$ block of $\mathsf X_\gamma-\lambda I$
leaves, in coordinates $(x,q,a,b)$,
\begin{equation}
\begin{pmatrix}
 \gamma^2S_{\gamma,\lambda}-\lambda I
    &\gamma C_2&\gamma^3C_1R&0\\
 -\gamma C_2^\top&(\gamma-\lambda)I&0&-\gamma^2I\\
 -\gamma^3RC_1^\top&0&(\gamma^3-\lambda)I-\gamma^4R&0\\
 0&-\gamma^2I&0&(\gamma^3-\lambda)I
\end{pmatrix}.
\label{eq:proof-first-elimination}
\end{equation}
Multiplying the determinant of this matrix by
$\det(D+\gamma I-\lambda I)$ gives
$\det(\mathsf X_\gamma-\lambda I)$.
Set
\[
    \mathsf Y=\begin{pmatrix}0&E\\-E^\top&I\end{pmatrix}.
\]
At $\lambda=\gamma\zeta$, dividing each row of
\eqref{eq:proof-first-elimination} by $\gamma$ gives
\begin{equation}
    \gamma^{-(d-r)}\det(\mathsf X_\gamma-\gamma\zeta I)
    \longrightarrow
    \det D\det(\mathsf Y-\zeta I)(-\zeta)^{m+d_2},
    \label{eq:proof-first-eigenvalue-scale}
\end{equation}
locally uniformly in complex $\zeta$.
An orthogonal singular-value decomposition of $E$ reduces $\mathsf Y$
to $k$ blocks $\bigl(\begin{smallmatrix}0&\sigma_i\\-\sigma_i&1\end{smallmatrix}\bigr)$.
Hence
\[
    \det(\mathsf Y-\mu I)=\prod_{i=1}^k(\mu^2-\mu+\sigma_i^2),
    \ \sigma_i=\sigma_i(E)=\sigma_i(C_2)>0.
\]
Each root has positive real part. Rouch\'e's theorem applied
to disjoint circles around the distinct roots in
\eqref{eq:proof-first-eigenvalue-scale} gives $2k$ eigenvalues with $\lambda/\gamma$ tending to
these roots, counted with multiplicity.

Uniformly for complex $|\lambda|\le c\gamma^2$, where $c$ is fixed,
the $(\xi,q)$ block of \eqref{eq:proof-first-elimination}
has inverse
\[
    (\gamma\mathsf Y+O(\gamma^2))^{-1}
    =\gamma^{-1}
      \begin{pmatrix}E^{-\top}E^{-1}&-E^{-\top}\\E^{-1}&0\end{pmatrix}
       +O(1).
\]
Its $qq$ block is $O(1)$. The Schur correction on $(u,a,b)$ is
\[
\begin{aligned}
    &\begin{pmatrix}O(\gamma^2)&0\\O(\gamma^3)&0\\0&O(\gamma^2)\end{pmatrix}
     \begin{pmatrix}O(\gamma^{-1})&O(\gamma^{-1})\\O(\gamma^{-1})&O(1)\end{pmatrix}
     \begin{pmatrix}O(\gamma^2)&O(\gamma^3)&0\\0&0&O(\gamma^2)\end{pmatrix}\hspace{1em}=
      \begin{pmatrix}
       O(\gamma^3)&O(\gamma^4)&O(\gamma^3)\\
       O(\gamma^4)&O(\gamma^5)&O(\gamma^4)\\
       O(\gamma^3)&O(\gamma^4)&O(\gamma^4)
      \end{pmatrix}.
\end{aligned}
\]
Since $S_{\gamma,\lambda}=S+O(\gamma)$, subtracting this
correction and combining $(a,b)=v$ leaves
\begin{equation}
    \begin{pmatrix}
       \gamma^2S_0-\lambda I+O(\gamma^3)&O(\gamma^3)\\
       O(\gamma^3)&(\gamma^3-\lambda)I_{d_2}+O(\gamma^4)
    \end{pmatrix}.
    \label{eq:proof-second-elimination}
\end{equation}
Taking $\lambda=\gamma^2\zeta$ gives, locally uniformly,
\begin{equation}
    \gamma^{-2(k+m+d_2)}
    \det(\mathsf X_\gamma-\gamma^2\zeta I)
    \longrightarrow
    \det D\det\mathsf Y\det(S_0-\zeta I)(-\zeta)^{d_2}.
    \label{eq:proof-second-eigenvalue-scale}
\end{equation}
The nonzero roots corresponding to $S_0\succ0$ therefore
give $m$ eigenvalues whose ratios $\lambda/\gamma^2$ tend
to the positive eigenvalues of $S_0$, by Rouch\'e's theorem.

At $\lambda=\gamma^3\zeta$, the upper-left block in
\eqref{eq:proof-second-elimination} has inverse
$\gamma^{-2}(S_0^{-1}+O(\gamma))$.
Its elimination changes the lower-right block by
\[
    O(\gamma^3)\,O(\gamma^{-2})\,O(\gamma^3)=O(\gamma^4),
\]
so the final block is $\gamma^3((1-\zeta)I_{d_2}+O(\gamma))$.
Consequently,
\begin{equation}
    \gamma^{-(2k+2m+3d_2)}
    \det(\mathsf X_\gamma-\gamma^3\zeta I)
    \longrightarrow
    \det D\det\mathsf Y\det S_0(1-\zeta)^{d_2}.
    \label{eq:proof-third-eigenvalue-scale}
\end{equation}
Rouch\'e's theorem gives $d_2$ roots with
$\lambda/\gamma^3\to1$.
The four groups are disjoint for small $\gamma$ and exhaust
the spectrum because $r+2k+m+d_2=d$.
When a block is empty, its elimination and root group are omitted.

Since $\norm{H}\le1$, the limits $\delta_i$ lie in $(0,1]$,
and the first group satisfies
\[
    |1-\lambda_i(\gamma)/4|\longrightarrow1-\delta_i/4<1.
\]
For each other group, there are $a,b>0$ and $j\in\{1,2,3\}$
such that
\[
    \operatorname{Re}\lambda\ge a\gamma^j,\ 
    |\lambda|\le b\gamma^j,\ 
    |1-\lambda/4|^2
    \le1-\frac a2\gamma^j+\frac{b^2}{16}\gamma^{2j}<1
\]
for sufficiently small $\gamma$. This proves the claim.
\end{proof}

\begin{proof}[Proof of Lemma \ref{lem:fixed-parameter-convergence}]
The left inequality in the calm local minimax definition gives
$\nabla_y f(z^*)=0$ and $B\preceq0$.
For any unit $u$ and the response constant $M$,
\[
    0\le\max_{\norm{h}\le M\delta}
       \bigl(f(x^*+\delta u,y^*+h)-f(z^*)\bigr)
      =\delta\nabla_x f(z^*)^\top u+O(\delta^2).
\]
Applying this to $u$ and $-u$ gives $\nabla_x f(z^*)=0$.
Lemma~\ref{lem:curvature-characterization} yields $S_0\succ0$.

Differentiating the three gradient updates at $w^*$ gives
\[
    \mathsf J_\gamma=
    \begin{pmatrix}
       I-\gamma^2A/4&-\gamma^2C/4&0\\
       C^\top/4&I+B/4-\gamma I/4&\gamma I/4\\
       0&\gamma^3I/4&I-\gamma^3I/4
    \end{pmatrix}.
\]
With $D_\gamma=\operatorname{diag}(\gamma I,I,\gamma I)$, $D_\gamma^{-1}\mathsf J_\gamma D_\gamma
    =I-\mathsf X_\gamma/4$.

By Lemma~\ref{lem:spectral-stability}, every eigenvalue of
$\mathsf J_\gamma$ has modulus less than one for sufficiently small $\gamma$.
The gradient updates are $C^1$ and fix $w^*$.
Lemma~\ref{lem:tools-linear-stability} gives
$C_0\ge1$, $q\in(0,1)$, and $r_{\mathrm{stab}}>0$ such that
\[
    \norm{w_t-w^*}\le C_0q^t\norm{w_0-w^*}
    \ \text{if }\norm{w_0-w^*}<r_{\mathrm{stab}}
\]
for the gradient updates alone.

We finish by checking that the comparison step leaves these
updates unchanged near $w^*$. Local maximality in $y$ supplies
a radius $\rho>0$ such that
$f(x^*,y^*+h)\le f(z^*)$ whenever $\norm{h}\le\rho$.
Reduce the upper bound on $\gamma$ so that $\ell=\gamma<\rho$.
The derivative at $w^*$ of
\[
    (x,y,v)\longmapsto
       \bigl(\nabla f(x,y),s(y-v)\bigr)
\]
is nonsingular: its determinant is
$(-s)^{d_2}\det\nabla^2 f(z^*)$.
Lemma~\ref{lem:tools-error-bound} and equivalence of the
Euclidean norm and the sum of the two block norms give
$c_*>0$ and a neighborhood on which
\begin{equation}
    \norm{\nabla f(x,y)}+s\norm{y-v}
       \ge c_*\norm{w-w^*}.
    \label{eq:proof-residual-lower-bound}
\end{equation}

Let $e=\norm{w-w^*}$, and let
$(x',\widetilde y',v')$ be the result of one gradient update
from $w$. Since $s\le L$ and $\tau\ge1$,
\[
\begin{aligned}
    \norm{(x',\widetilde y',v')-w}
    &\le \eta\norm{\nabla f(x,y)}
       +\eta s\sqrt{1+\tau^{-2}}\norm{y-v}\le \eta(L+2s)e\le\tfrac34e.
\end{aligned}
\]
In particular,
$\norm{(x',\widetilde y')-z^*}\le2e$.
For any $\norm{h}\le\ell$, the function
$z\mapsto f(z+(0,h))-f(z)$ is $L\norm{h}$-Lipschitz,
by Assumption~\ref{ass:smooth}.
Its value at $z^*$ is nonpositive, so
\[
    f(x',\widetilde y'+h)-f(x',\widetilde y')
       \le 2L\ell e.
\]
At the current state, \eqref{eq:proof-residual-lower-bound}
and \eqref{eq:adaptive-margin} bound the comparison margin:
\[
    \vartheta\ge\ell\sqrt{L\ell c_*e}.
\]
Thus every comparison returns zero whenever
$e\le c_*\ell/(4L)$ and $w$ lies in the neighborhood where
\eqref{eq:proof-residual-lower-bound} holds.
At $w=w^*$ the same conclusion follows directly from local
maximality and the strict inequality in the comparison oracle.

Choose $r_{\mathrm{rej}}>0$ so that the comparison conclusion
holds whenever $\norm{w-w^*}<r_{\mathrm{rej}}$, and take
\[
    r_{\mathrm{init}}=
       \min\{r_{\mathrm{stab}},r_{\mathrm{rej}}/C_0\}.
\]
For $\norm{w_0-w^*}<r_{\mathrm{init}}$, every gradient iterate
satisfies
\[
    \norm{w_t-w^*}
       \le C_0q^t\norm{w_0-w^*}
       <C_0r_{\mathrm{init}}\le r_{\mathrm{rej}}.
\]
Induction on $t$ therefore shows that every comparison is
rejected and the complete iterates coincide with the gradient
iterates, for every sequence $(U_t)$ satisfying $\norm{U_t}\le1$.
This proves the lemma with initial radius $r_{\mathrm{init}}$.
\end{proof}

\begin{proof}[Proof of Lemma \ref{lem:continuing-trial}]
Fix a trial $j$ that continues forever and write
$\alpha=\alpha_j$, $s=s_j$, and $\tau=\tau_j$.
Since the trial starts only when $R_j>0$, its parameter
satisfies $0<\alpha\le1/2$.
For integers $m>n\ge0$, the restart check gives
\begin{align*}
    \norm{w_{t_j+m}-w_{t_j+n}}
    &\le
    \sum_{k=n}^{m-1}
        \norm{w_{t_j+k+1}-w_{t_j+k}}\le
    \sum_{k=n}^{m-1}\alpha^2(1-\alpha)^k\le \alpha(1-\alpha)^n.
\end{align*}
Consequently, the sequence is Cauchy. Denote its limit by
$(x^\infty,y^\infty,v^\infty)$ and let $m\to\infty$
in this inequality to obtain
\begin{equation}
    \norm{w_{t_j+n}-(x^\infty,y^\infty,v^\infty)}
    \le \alpha(1-\alpha)^n.
    \label{eq:continuing-trial-tail}
\end{equation}

The comparison step does not change the proposed $x$ and $v$
coordinates. Thus, \eqref{eq:base-update} gives
\begin{align*}
    \nabla_x f(x_t,y_t)
    =\frac{\tau}{\eta}(x_t-x_{t+1})
      \longrightarrow0,\ 
    y_t-v_t
    =\frac{\tau}{\eta s}(v_{t+1}-v_t)
      \longrightarrow0.
\end{align*}
By continuity,
\[
    \nabla_x f(x^\infty,y^\infty)=0,
    \  v^\infty=y^\infty.
\]
It remains to prove that the other gradient component vanishes.
Set $g=\nabla_y f(x^\infty,y^\infty)$.
The provisional maximizer update satisfies
\[
    \widetilde y_{t+1}
    =y_t+\eta\bigl(\nabla_y f(x_t,y_t)-s(y_t-v_t)\bigr)
    \longrightarrow y^\infty+\eta g.
\]
Suppose $g\ne0$. If infinitely many comparison queries returned
zero, then, along those indices,
\[
    y_{t+1}-y_t
    =\eta\bigl(\nabla_y f(x_t,y_t)-s(y_t-v_t)\bigr)
    \longrightarrow\eta g\ne0,
\]
contradicting convergence of $y_t$.
Therefore, all sufficiently late queries return one. Their
acceptance inequalities and $\vartheta_t\ge0$ imply
\[
    f(x_{t+1},y_{t+1})
    >f(x_{t+1},\widetilde y_{t+1})+\vartheta_t
    \ge f(x_{t+1},\widetilde y_{t+1}).
\]
Passing to the limit yields
\begin{equation}
    f(x^\infty,y^\infty)
    \ge f(x^\infty,y^\infty+\eta g).
    \label{eq:continuing-trial-comparison-limit}
\end{equation}
On the other hand, the Hessian bound in
Assumption~\ref{ass:smooth} gives
\begin{align*}
    f(x^\infty,y^\infty+\eta g)-f(x^\infty,y^\infty)
    &\ge \eta\norm{g}^2
        -\frac{L\eta^2}{2}\norm{g}^2=\eta\left(1-\frac{L\eta}{2}\right)\norm{g}^2
      >0,
\end{align*}
where $L\eta=1/4$.
This contradicts \eqref{eq:continuing-trial-comparison-limit}.
Hence $g=0$, and the limit in
\eqref{eq:continuing-trial-tail} is the stationary lift
$w^\infty=(x^\infty,y^\infty,y^\infty)$.
\end{proof}

\begin{proof}[Proof of Lemma \ref{lem:successful-trial}]
Choose $\gamma_*$ as in
Lemma~\ref{lem:fixed-parameter-convergence}, and fix
$\gamma\in(0,\gamma_*)$.
Throughout this proof,
\[
    s=L\gamma,\  \tau=\gamma^{-2},\ 
    \ell=\gamma,\  \eta=\frac1{4L}.
\]
That lemma supplies $r_0>0$, $C\ge1$, and $q\in(0,1)$
such that the updates without restarts satisfy
\begin{equation}
    \norm{w_{t_j+n}-w^*}
    \le Cq^n\norm{w_{t_j}-w^*},
    \  n\ge0,
    \label{eq:successful-trial-fixed-rate}
\end{equation}
whenever $\norm{w_{t_j}-w^*}<r_0$.
Every comparison then returns zero for every sequence $(U_t)$
satisfying $\norm{U_t}\le1$.

We next relate the initialization perturbation to its center's
distance from $w^*$.
The derivative at $w^*$ of the map $    (x,y,v)\longmapsto
    \begin{pmatrix}\nabla f(x,y)\\s(y-v)\end{pmatrix}$ is $
    \begin{pmatrix}
        \nabla^2 f(z^*) & 0\\
        \begin{matrix}0&sI_{d_2}\end{matrix} & -sI_{d_2}
    \end{pmatrix}$.

Its determinant is
$(-s)^{d_2}\det\nabla^2 f(z^*)\ne0$.
Lemma~\ref{lem:tools-error-bound}, together with
$\sqrt{a^2+b^2}\le a+b\le\sqrt2\sqrt{a^2+b^2}$ for
$a,b\ge0$, therefore gives constants $c_s,C_s,r_c>0$
such that
\begin{equation}
    c_s\norm{w-w^*}
    \le \norm{\nabla f(x,y)}+s\norm{y-v}
    \le C_s\norm{w-w^*}
    \ \text{if }\norm{w-w^*}<r_c.
    \label{eq:successful-trial-residual}
\end{equation}

Define the constants
\[
    \kappa=1+\frac{L}{c_s},\ 
    a_0=\min\left\{
        \frac12,
        \left(\frac{r_0}{2\kappa}\right)^{1/3},
        \frac{1}{2C(1+q)\kappa},
        \frac{1-q}{2}
    \right\},
\]
and choose
\[
    0<\delta<\min\left\{r_c,\frac{La_0^3}{C_s}\right\}.
\]
Consider a center $\bar w_j$ with
$\norm{\bar w_j-w^*}<\delta$.
If $R_j=0$, the lower bound in
\eqref{eq:successful-trial-residual} implies
$\bar w_j=w^*$, so the algorithm returns $z^*$.
Otherwise, its restart parameter satisfies
\[
    0<R_j<La_0^3\le\frac{L}{8},
    \  \alpha_j^3=\frac{R_j}{L},
    \  0<\alpha_j<a_0.
\]
For every initialization draw with $\norm{Z_j}\le1$,
\begin{align*}
    \norm{w_{t_j}-w^*}
    &\le\norm{\bar w_j-w^*}+\alpha_j^3\norm{Z_j}\le\frac{R_j}{c_s}+\alpha_j^3
      =\kappa\alpha_j^3
      <\frac{r_0}{2}.
\end{align*}
We can therefore apply
\eqref{eq:successful-trial-fixed-rate} to the updates
without the restart test. For every $n\ge0$, they satisfy
\begin{align*}
    \norm{w_{t_j+n+1}-w_{t_j+n}}
    &\le\norm{w_{t_j+n+1}-w^*}
       +\norm{w_{t_j+n}-w^*}\le C(1+q)q^n\norm{w_{t_j}-w^*}\\
    &\le C(1+q)\kappa\alpha_j^3q^n\le\alpha_j^2(1-\alpha_j)^n.
\end{align*}
The last inequality uses
$C(1+q)\kappa\alpha_j\le1$ and $q\le1-\alpha_j$,
both of which follow from the choice of $a_0$.
Thus every proposed update passes
\eqref{eq:restart-check}. Induction on $n$ shows that
the actual trial follows these updates for all time.
Its limit is $w^*$ by
\eqref{eq:successful-trial-fixed-rate}, and no restart occurs.
\end{proof}

\begin{proof}[Proof of Lemma \ref{lem:local-restarts}]
Let $w^*=(x^*,y^*,y^*)$ be the stationary lift in the
statement and set $c_0=1+\sqrt2$.
For a reached trial $j\le K$, write
$e_j=\norm{\bar w_j-w^*}$.
Stationarity, the Hessian bound, and $s_j\le L$ give
\begin{align*}
    R_j
    &=\norm{\nabla f(\bar x_j,\bar y_j)}
       +s_j\norm{\bar y_j-\bar v_j}\\
    &\le L\norm{(\bar x_j,\bar y_j)-(x^*,y^*)}
       +L\norm{(\bar y_j-y^*)-(\bar v_j-y^*)}\\
    &\le Le_j+L\sqrt2\,e_j=c_0Le_j.
\end{align*}
If $R_j=0$, the algorithm returns at this center and
does not produce another trial. If $R_j>0$, then
\[
    \alpha_j^3
    =\min\left\{\frac18,\frac{R_j}{L}\right\}
    \le c_0e_j,
    \  \alpha_j\le(c_0e_j)^{1/3}.
\]
For every accepted state in this trial, including its
initialization when $n=0$, we have
\begin{align*}
    \norm{w_{t_j+n}-w^*}
    &\le e_j+\alpha_j^3\norm{Z_j}
        +\sum_{k=0}^{n-1}
            \norm{w_{t_j+k+1}-w_{t_j+k}}\\
    &\le e_j+\alpha_j^3
        +\sum_{k=0}^{n-1}\alpha_j^2(1-\alpha_j)^k\\
    &\le e_j+\alpha_j^3+\alpha_j\\
    &\le(1+c_0)e_j+(c_0e_j)^{1/3}.
\end{align*}
If this trial fails and $j<K$, the next center is the
retained state at failure. The same bound therefore gives
\begin{equation}
    e_{j+1}\le(1+c_0)e_j+(c_0e_j)^{1/3}.
    \label{eq:local-restart-displacement}
\end{equation}
In particular, these bounds do not depend on how long
the trial lasts or on the sampled vectors $U_t$.

To choose an initial neighborhood explicitly, define
\[
    \rho_0=\rho,
    \ 
    \rho_{m+1}
    =\min\left\{
        \frac{\rho_m}{2(1+c_0)},
        \frac{\rho_m^3}{8c_0}
    \right\},
    \  m=0,\ldots,J-1,
\]
and set $r=\rho_J>0$.
This sequence is decreasing, and, for $0\le e<\rho_{m+1}$,
\[
    (1+c_0)e+(c_0e)^{1/3}
    <\frac{\rho_m}{2}+\frac{\rho_m}{2}
    =\rho_m.
\]
Since $\bar w_1=w_0$, the hypothesis gives
$e_1<\rho_J$.
Induction using \eqref{eq:local-restart-displacement}
now yields, for each reached $j\le\min\{J,K\}$,
\[
    e_j<\rho_{J-j+1}\le\rho,
    \ 
    \norm{w_{t_j+n}-w^*}<\rho_{J-j}\le\rho
\]
for every state produced within that trial.
These are the required bounds for the centers,
initialization perturbations, and accepted iterates.
\end{proof}

\begin{proof}[Proof of Lemma \ref{lem:absolute-continuity}]
Write $\varepsilon=1/\tau$ and denote the displacement in
\eqref{eq:base-update} by
\[
    d(w)=\eta\begin{pmatrix}
        -\varepsilon\nabla_x f(x,y)\\
        \nabla_y f(x,y)-s(y-v)\\
        \varepsilon s(y-v)
    \end{pmatrix}.
\]
The gradient and anchor contributions give
\[
    \sup_w\norm{Dd(w)}
    \le\eta L+\eta s\sqrt{2(1+\varepsilon^2)}
    \le\eta(L+2s)\le\frac34.
\]
Consequently,
\begin{align*}
    \norm{w+d(w)-w'-d(w')}
    &\ge\norm{w-w'}-\norm{d(w)-d(w')}\ge\tfrac14\norm{w-w'},\\
    \norm{(I+Dd(w))h}&\ge\tfrac14\norm{h}.
\end{align*}
For each output $a$, the equation $w=a-d(w)$ has a unique
solution by contraction. The inverse function theorem then
makes $w\mapsto w+d(w)$ a global $C^1$ diffeomorphism,
with a $4$-Lipschitz inverse.

For a fixed realization $U_t=u$ and a null set $N$,
the comparison gives
\[
\begin{aligned}
    \{w:\text{the updated state belongs to }N\}\subseteq
    \{w:w+d(w)\in N\}
    \cup\{w:w+d(w)\in N-(0,\ell u,0)\}.
\end{aligned}
\]
Both sets on the right are null, since the inverse is
Lipschitz. Independence of $U_t$ and $w_t$ therefore yields
\[
    \Pr(w_{t+1}\in N)
    =\int\Pr(w_{t+1}\in N\mid U_t=u)\,d\Pr_{U_t}(u)=0
\]
whenever $w_t$ has an absolutely continuous law.
Induction from $w_0$ gives
\begin{equation}
    \Pr(w_t\in N)=0
    \ \text{for every null set $N$ and every finite $t$.}
    \label{eq:finite-iterate-absolute-continuity}
\end{equation}
\end{proof}

\begin{proof}[Proof of Lemma \ref{lem:gradient-instability}]
By Lemma~\ref{lem:curvature-characterization}, $S_0$ is nonsingular
and is not positive definite. Hence some
$u\in\ker C_2^\top$ satisfies $u^\top Su<0$.
Put $\varepsilon=1/\tau$ and $\bar s=s/L$.
For $q>0$, define
\[
    \beta(q)=q+\frac{\bar s q}{\varepsilon\bar s+q},
    \  Q(q)=A+\frac q\varepsilon I
        +C(-B+\beta(q)I)^{-1}C^\top.
\]
The inverse exists because $-B\succeq0$ and $\beta(q)>0$.
Since $C_2^\top u=0$,
\begin{align*}
    u^\top Q(q)u
    &=u^\top Au+\frac q\varepsilon\norm{u}^2
      +u^\top C_1(D+\beta(q)I)^{-1}C_1^\top u\longrightarrow u^\top Su<0.
\end{align*}
For $q>\varepsilon\norm{A}$, $Q(q)\succ0$.
Continuity of its smallest eigenvalue therefore gives
$q_*>0$ and $a\ne0$ such that $Q(q_*)a=0$.
Set
\[
    b=(-B+\beta(q_*)I)^{-1}C^\top a,
    \  c=\frac{\varepsilon\bar s}{\varepsilon\bar s+q_*}b.
\]
The Jacobian is
\[
    \mathsf J=I-\frac14
    \begin{pmatrix}
        \varepsilon A&\varepsilon C&0\\
        -C^\top&\bar sI-B&-\bar sI\\
        0&-\varepsilon\bar sI&\varepsilon\bar sI
    \end{pmatrix}.
\]
The definitions imply
\begin{align*}
    \varepsilon Aa+\varepsilon Cb=-q_*a,\ 
    -C^\top a+(\bar sI-B)b-\bar sc=-q_*b,\ 
    -\varepsilon\bar sb+\varepsilon\bar sc=-q_*c.
\end{align*}
Therefore
\[
    \mathsf J\begin{pmatrix}a\\b\\c\end{pmatrix}
    =\left(1+\frac{q_*}{4}\right)
      \begin{pmatrix}a\\b\\c\end{pmatrix},
    \ 1+\frac{q_*}{4}>1.
\]
\end{proof}

\begin{proof}[Proof of Lemma \ref{lem:avoidance-fixed-parameters}]
Lemma~\ref{lem:absolute-continuity} gives absolute continuity
of every finite iterate. Fix a nondegenerate stationary point
$z^*=(x^*,y^*)$ that is not calm local minimax, and set
$w^*=(x^*,y^*,y^*)$.

Suppose first that some $h$ satisfies
\[
    0<\norm{h}<\ell,\ 
    \Delta=f(x^*,y^*+h)-f(z^*)>0.
\]
Continuity supplies an open ball $\mathcal U$ around $h/\ell$,
with closure inside the unit ball, such that
\[
    f(x^*,y^*+\ell u)-f(z^*)>3\Delta/4,
    \  \ell\norm{u}>3\norm{h}/4,
    \  u\in\overline{\mathcal U}.
\]
The gradient candidate tends to $w^*$ and the margin tends
to zero as $w_t\to w^*$.
Hence some $0<\rho<\norm{h}/2$ satisfies the following:
if $\norm{w_t-w^*}<\rho$, then, uniformly over
$u\in\overline{\mathcal U}$, we have
\begin{align*}
    f(x_{t+1}^{\mathrm{cand}},\widetilde y_{t+1}+\ell u)
      -f(x_{t+1}^{\mathrm{cand}},\widetilde y_{t+1})
      >\Delta/2>\vartheta_t,\ 
    \norm{\widetilde y_{t+1}-y^*}&<\norm{h}/4.
\end{align*}
The corresponding candidate action is accepted and gives
\[
    \norm{y_{t+1}-y^*}
    \ge\ell\norm{u}-\norm{\widetilde y_{t+1}-y^*}
    >\norm{h}/2>\rho.
\]
Writing $p=\Pr(U_t\in\mathcal U)>0$, independence of
$U_t$ from the preceding history implies
\[
    \Pr\!\left(\norm{w_t-w^*}<\rho
        \text{ for }t=n,\ldots,n+m\right)\le(1-p)^m.
\]
Therefore
\[
    \Pr(w_t\to w^*)
    \le\sum_{n\ge0}\Pr\!\left(
       \norm{w_t-w^*}<\rho\text{ for all }t\ge n
    \right)=0.
\]

Otherwise, continuity gives
\begin{equation}
    f(x^*,y^*+h)\le f(z^*),\ \norm{h}\le\ell.
    \label{eq:avoidance-inner-maximum}
\end{equation}
The derivative at $w^*$ of
$(x,y,v)\mapsto(\nabla f(x,y),s(y-v))$ is
\[
    \begin{pmatrix}
        \nabla^2f(z^*)&0\\
        s\begin{pmatrix}0&I\end{pmatrix}&-sI
    \end{pmatrix},
\]
whose determinant is $(-s)^{d_2}\det\nabla^2f(z^*)\ne0$.
Lemma~\ref{lem:tools-error-bound} supplies $c_*>0$ such that,
for $e=\norm{w_t-w^*}$ sufficiently small,
\[
    \norm{\nabla f(x_t,y_t)}+s\norm{y_t-v_t}\ge c_*e,
    \  \vartheta_t\ge\ell\sqrt{L\ell c_*e}.
\]
The gradient update is $7/4$-Lipschitz by the proof of
Lemma~\ref{lem:absolute-continuity} and fixes $w^*$.
Moreover, $z\mapsto f(z+(0,h))-f(z)$ is $L\norm{h}$-Lipschitz.
Together with \eqref{eq:avoidance-inner-maximum}, this yields
\begin{align*}
    f(x_{t+1}^{\mathrm{cand}},\widehat y_{t+1})
      -f(x_{t+1}^{\mathrm{cand}},\widetilde y_{t+1})\hspace{1em}\le L\ell\norm{(x_{t+1}^{\mathrm{cand}},\widetilde y_{t+1})-z^*}
      \le2L\ell e\le\ell\sqrt{L\ell c_*e}\le\vartheta_t
\end{align*}
whenever $0<e<c_*\ell/(4L)$ and the residual bound holds.
At $e=0$ the gain is nonpositive and the margin is zero.
Thus every comparison returns zero near $w^*$.

Lemma~\ref{lem:gradient-instability} gives an eigenvalue
greater than one for the Jacobian of the gradient update.
This update is a $C^1$ diffeomorphism by the proof of
Lemma~\ref{lem:absolute-continuity}.
Lemma~\ref{lem:tools-center-stable} consequently gives a
neighborhood $V$ and an embedded manifold $\mathcal W$ such that
\[
    \dim\mathcal W<d_1+2d_2,\ 
    \operatorname{Leb}(\mathcal W)=0,
\]
and every gradient trajectory remaining in $V$ starts in
$\mathcal W$.
Shrink $V$ inside the neighborhood where every comparison
returns zero. Convergence to $w^*$ implies that, after some
finite step, the iterates remain in $V$ and follow the gradient
updates. Hence, by Lemma~\ref{lem:absolute-continuity}, we have
\[
    \{w_t\to w^*\}\subseteq\bigcup_{n\ge0}\{w_n\in\mathcal W\},
    \ 
    \Pr(w_t\to w^*)\le\sum_{n\ge0}\Pr(w_n\in\mathcal W)=0
\]
\end{proof}

\begin{proof}[Proof of Lemma \ref{lem:correct-trial-limits}]
We first show that the stationarity checks cannot stop the
algorithm with positive probability.
For any $s>0$,
\begin{equation}
    \norm{\nabla f(x,y)}+s\norm{y-v}=0
    \ \longrightarrow\  y=v.
    \label{eq:zero-residual-null-set}
\end{equation}
The set $\{(x,y,v):y=v\}$ is a linear subspace of
$\RR^{d_1+2d_2}$ of codimension $d_2\ge1$, and therefore
has Lebesgue measure zero.
The first center is $w_0=(x_0,y_0,y_0)$, for which
$R_1=\norm{\nabla f(x_0,y_0)}>0$ by assumption.

Consider any reached trial $j$ with $R_j>0$, conditional
on its center, parameters, and the history before drawing $Z_j$.
These quantities determine $\alpha_j>0$.
The initialization
\[
    w_{t_j}=\bar w_j+\alpha_j^3Z_j
\]
is uniform on a ball of positive radius, so its conditional
law is absolutely continuous.
For this argument only, continue the updates of this trial
without restarting even if its check fails, and denote
these states by $u_0,u_1,\ldots$, with $u_0=w_{t_j}$.
Lemma~\ref{lem:avoidance-fixed-parameters} gives
\begin{equation}
    \Pr(u_n\in N\mid\text{history before }Z_j)=0,
    \  n\ge0,
    \label{eq:conditional-null-states}
\end{equation}
for every Lebesgue-null set $N$.
The random vectors $U_t$ used in this continuation are
independent of the initialization and have the same law
as those used by the algorithm.
Until a check fails, $u_n$ coincides with the actual state
at step $t_j+n$.

Let $\sigma$ be the first index at which this trial fails
its check, with $\sigma=\infty$ if it never fails.
On $\{\sigma<\infty\}$ the discarded candidate is
$u_{\sigma+1}$, and the retained state is $u_\sigma$.
Thus \eqref{eq:conditional-null-states} implies
\begin{align*}
    \Pr(u_\sigma\in N,\ \sigma<\infty
         \mid\text{history before }Z_j)&
    =\sum_{n=0}^{\infty}
      \Pr(u_n\in N,\ \sigma=n
          \mid\text{history before }Z_j)\\
    &\hspace{1em}
    \le\sum_{n=0}^{\infty}
      \Pr(u_n\in N\mid\text{history before }Z_j)=0.
\end{align*}
If $j+1\le K$, this retained state is the next center.
If $j+1>K$, the next center has the Gaussian density
specified in \eqref{eq:global-trial-sampling}.
In either case, \eqref{eq:zero-residual-null-set} gives
\[
    \Pr(\text{trial }j+1\text{ is reached and }R_{j+1}=0)=0.
\]
Induction from $R_1>0$ and a countable union over $j$ yield
\begin{equation}
    \Pr(R_j>0\text{ for every reached trial }j)=1.
    \label{eq:positive-residual-all-trials}
\end{equation}

We next identify the possible limits of a trial that continues
forever. At every stationary point $z^*$, the derivative of
$\nabla f$ is the nonsingular matrix $\nabla^2 f(z^*)$.
Lemma~\ref{lem:tools-error-bound} gives a neighborhood
containing no other stationary point.
Using a countable basis of open balls in $\RR^{d_1+d_2}$,
assign to each stationary point a basis ball containing it
and no other stationary point.
Different stationary points are assigned different balls;
therefore the set of stationary points is countable.
In particular, its subset of points that are not calm local
minimax is countable; denote this subset by $\mathcal B$.

Fix a reached trial $j$ and condition on its history before
$Z_j$, with $R_j>0$ as established in
\eqref{eq:positive-residual-all-trials}.
For each $z^*=(x^*,y^*)\in\mathcal B$,
Lemma~\ref{lem:avoidance-fixed-parameters} applied to
the continuation $u_n$ gives
\[
    \Pr\!\left(
        \lim_{n\to\infty}u_n=(x^*,y^*,y^*)
        \,\middle|\,\text{history before }Z_j
    \right)=0.
\]
It follows that
\begin{align*}
    &\Pr\!\left(
        \sigma=\infty,\ 
        \lim_{n\to\infty}u_n
        \in\{(x,y,y):(x,y)\in\mathcal B\}
        \,\middle|\,\text{history before }Z_j
    \right)\\
    \le&\sum_{(x^*,y^*)\in\mathcal B}
      \Pr\!\left(
        \lim_{n\to\infty}u_n=(x^*,y^*,y^*)
        \,\middle|\,\text{history before }Z_j
      \right)=0.
\end{align*}
On $\{\sigma=\infty\}$,
Lemma~\ref{lem:continuing-trial} already ensures that
$u_n$ converges to a stationary lift.
The preceding inequality excludes every incorrect
stationary lift.
Integrating over the finite history and taking a countable
union over the trial index $j$ proves the result.
\end{proof}

\begin{proof}[Proof of Lemma \ref{lem:global-restarts}]
Let $z^*=(x^*,y^*)$ be a nondegenerate calm local minimax
point, and set $w^*=(x^*,y^*,y^*)$.
By Lemma~\ref{lem:successful-trial}, there are an integer
$m_*\ge1$ and a radius $\delta>0$ such that a trial
with scale $2^{-m_*}$ and center within distance $\delta$
of $w^*$ never restarts.
For this fixed scale, the same lemma identifies a zero
residual in that ball only at $w^*$.
Consequently, every center in the nonempty open set
\[
    V=\{\bar w:\delta/4<\norm{\bar w-w^*}<\delta/2\}
\]
starts a trial that continues forever and converges to $w^*$.

Write $d=d_1+2d_2$.
At every globally sampled trial, independence of the center
and scale gives the same probability of sampling a center
and scale with the preceding property:
\begin{equation}
    p
    =2^{-m_*}\int_V
      \frac{1}{(2\pi)^{d/2}}
      \exp\!\left(-\frac{\norm{w-w_0}^2}{2}\right)\,dw
    >0.
    \label{eq:global-success-probability}
\end{equation}
The integral is positive because its integrand is positive
everywhere and $V$ is a nonempty open set.
The fixed Gaussian center $w_0$ makes $p$ independent of
the number and outcomes of preceding trials.

To bound the number of failures, sample in advance the
independent center and scale for each possible trial
$K+h$, $h\ge1$; unused samples do not affect the algorithm.
Let $A_h$ be the event that this center belongs to $V$
and this scale equals $2^{-m_*}$.
Then
\[
    \Pr(A_h)=p,\ 
    \Pr\!\left(\bigcap_{h=1}^{n}A_h^c\right)=(1-p)^n.
\]
If trial $K+h$ is reached and $A_h$ occurs, that trial
does not restart. Therefore
\[
    \{\text{trial }K+n+1\text{ is reached}\}
    \subseteq\bigcap_{h=1}^{n}A_h^c
\]
and
\[
    \Pr(\text{trial }K+n+1\text{ is reached})
    \le(1-p)^n\longrightarrow0.
\]
Since $\nabla f(x_0,y_0)\ne0$, the initial value $K$
is finite. An infinite number of restarts would reach
every trial $K+n+1$, so its probability is zero.
After finitely many restarts, the algorithm thus either
terminates at a stationarity check or continues forever
in its current trial.
\end{proof}

\end{document}